\documentclass{article}
\usepackage{array}
\usepackage{amsmath}
\usepackage{amssymb}
\usepackage{amsthm}
\usepackage{amsfonts}
\usepackage{graphicx}
\usepackage[toc,page]{appendix}
\usepackage[margin=0pt,font+=small,labelformat=parens,labelsep=space, skip=6pt,list=false,hypcap=false] {caption}
\usepackage{subfig}
\usepackage{amsthm}
\newtheorem{theorem}{Theorem}[section]
\usepackage[top=.9in,bottom=.9in]{geometry}
 \numberwithin{equation}{section}
 
\begin{document}

%
\let\labeldefs\iffalse
\pagestyle{plain}
\numberwithin{equation}{section}
\begin{flushleft} \vskip 0.3 in 
\iftrue
\centerline{A Series Representation for Riemann's Zeta Function and some Interesting Identities that Follow.} \vskip .3in 
\vskip .2in
\centerline{ Michael Milgram\footnote{mike@geometrics-unlimited.com}}
\centerline{Consulting Physicist, Geometrics Unlimited, Ltd.}
\centerline{Box 1484, Deep River, Ont. Canada. K0J 1P0}
\centerline{}
\centerline{ Sept. 9, 2020}
\centerline{}
\centerline{}
\vskip .1in
\fi
\centerline{}
\vskip .1in
MSC classes: 	11M06, 11M26, 11M35, 11M99, 26A09, 30B40, 30E20, 33C20, 33B20, 33B99
\vskip 0.1in 

\centerline{\bf Abstract}\vskip .3in

Using Cauchy's Integral Theorem as a basis, what may be a new series representation for Dirichlet's function $\eta(s)$, and hence Riemann's function $\zeta(s)$, is obtained in terms of the Exponential Integral function $E_{s}(i\kappa)$ of complex argument. From this basis, infinite sums are evaluated, unusual integrals are reduced to known functions and interesting identities are unearthed. The incomplete functions $\zeta^{\pm}(s)$ and $\eta^{\pm}(s)$ are defined and shown to be intimately related to some of these interesting integrals. An identity relating Euler, Bernouli and Harmonic numbers is developed. It is demonstrated that a known simple integral with complex endpoints can be utilized to evaluate a large number of different integrals, by choosing varying paths between the endpoints.


\section{Introduction}
In a recent paper \cite{Milgram2019}, I have developed an integral equation for Riemann's function $\xi(s)$, based on LeClair's series \cite{LeClair:2013oda} involving the Generalized Exponential Integral. Glasser, in a subsequent paper \cite{Glasser2019note}, has shown that this integral equation is equivalent to the standard statement of Cauchy's Integral Theorem applied to $\xi(s)$, from which he closed the circle \cite[Eq. (13)]{Glasser2019note} by re-deriving LeClair's original series representation. This raises the question of what will be found if Cauchy's Integral Theorem is applied to other functions related to Riemann's function $\zeta(s)$, notably the so-called Dirichlet function $\eta(s)$ (sometimes called the alternating zeta function), working backwards to a series representation? That question is the motivation  for this work.\newline

Section \ref{sec:Series} gives some preliminaries and develops the reduction of Cauchy's Integral Theorem to yield a series representation of $\eta(s)$ in terms of (generalized) Integro-Exponential functions \cite{Milgram:1985} $E_s(\kappa)$ (a.k.a Incomplete Gamma functions - see \eqref{EiGamma}). Although it is well-known that $\zeta(s)$ can be expressed \cite{LeClair:2013oda}, \cite{ParisExp} in terms of such functions in various ways, (see Appendix \ref{sec:EiSums}) the form derived here differs substantially from these since it involves a series that indexes half-integer values of the imaginary argument of $E_s(i\kappa)$ and does not possess generic $s\leftrightarrow (1-s)$ symmetry. Several subsections explore variations of this series, revealing some interesting identities, not all of which pertain to either $\zeta(s)$ or $\eta(s)$. Section \ref{sec:Derivative} presents a brief analysis of these results as they pertain to $\zeta^{\prime}(s)$.  In Section \ref{sec:Integro}, the series representation is converted into an integral representation whose endpoints are complex, and in a number of subsections this representation is reduced to a real integral representation by choosing a number of different paths between the endpoints in the complex plane, yielding an evaluation of several unusual integrals. Appendix \eqref{sec:Def_and_Lemmas} collects some known results referenced in the main text, Appendix \eqref{sec:Theorem} contains the proof of a needed theorem and Appendix \eqref{sec:Circuitous} presents the details of a lengthy derivation required in Subsection \eqref{sec:s_2mP1_Lt}. \newline

In the late 1800's, Russell (\cite{RussellA}) initiated a tradition by presenting a list of 12 integral identities with little in the way of proof. This tradition has been recently (2008) reinvigorated by Amdeberhan and Moll(\cite{AmMoll2008dozen}); more recently (2018) Coffey (\cite{Coffey_2018}) has provided some proofs and generalizations of Russell's identities. Continuing the tradition, various subsections throughout this work present (special case) integrals that satisfy Russell's presumed criterion - that being that they have no apparent application other than that they are interesting and involve complicated trigonometric and hyperbolic functions. See also \cite{Coffey2018integrals}.\newline

Therefore, a (baker's) dozen of the most appealing of the interesting results obtained herein follows immediately below, cross-referenced to where each appears in the text of this paper. Departing from Russell's approach, each cross-reference to the later text, provides sufficient information to allow an informed reader to obtain a formal proof if (s)he so desires. A list of notations appears in 
Appendix \ref{sec:AppendixA}. \newline
\begin{itemize}

\item 
The main result upon which this paper is based (see \eqref{R5b}):
\begin{align} \nonumber
\displaystyle & \left( {2}^{s}-1 \right) \zeta \left( s \right) =-{\pi}^{s-1}\sin \left( \pi s/2 \right)\sum _{j=0}^{\left\lfloor (n-1)/2\right\rfloor}{\frac {\Gamma \left( 2\,j+1-s \right) 
\mbox{}{\mathcal{E}(2j)} }{\Gamma \left( 2\,j+1 \right) }} \\ \nonumber
&+{\frac {{\pi}^{s-n}{{\rm e}^{i\pi\,n/2}}}{2\Gamma \left( s-n \right) \cos \left( \pi s/2 \right) }\sum _{k=0}^{\infty }{\frac {{E_{1+n-s}} \left(-i\pi\, \left( k+1/2 \right)  \right) + \left( -1 \right) ^{n}{E_{1+n-s}} \left(i\pi\, \left( k+1/2 \right)  \right) 
\mbox{}}{ \left( k+1/2 \right) ^{n}}}}
\end{align} 

\item 
An unusual recursion relation for Bernoulli numbers (see \eqref{BernId}):

\begin{equation} \nonumber
\displaystyle {B_{2N+2}}  ={\frac {\Gamma \left( 3+2\,N \right) }{{2}^{2\,N} \left( {2}^{2\,N+2}-1 \right) }\sum _{k=0}^{N-1}{\frac {{2}^{2\,k} \left( 1-{2}^{2\,k+2} \right) B_{2k+2} }{\Gamma \left( 1-2\,k+2\,N \right) \Gamma \left( 2\,k+3 \right) 
\mbox{}}}}+{\frac {1+N}{{2}^{2\,N+1} \left( {2}^{2\,N+2}-1 \right) }}, \hspace{10pt} N\geq0\,, 
\end{equation} 

\item
The sum of an infinite series (see \eqref{T0aA}):
\begin{align} \nonumber
\displaystyle \frac{{\pi}^{p-2\,m} \left( -1 \right) ^{m}{{\rm e}^{-i\pi p/2}}}{\Gamma \left( 1+p \right)}
\mbox{}\sum _{k=0}^{\infty }{\frac {{ \left( -1 \right) ^{p}\,E_{-p}} \left( i\pi\, \left( k+1/2 \right)  \right) 
\mbox{}+{E_{-p}} \left(-i\pi\, \left( k+1/2 \right)  \right) }{ \left(k+1/2 \right) ^{-p+2\,m}}}\\=- \sum _{j=\left\lfloor -p/2+m+1/2\right\rfloor}^{m}{\frac {\mathcal{E}
 \left( 2\,j \right) }{\Gamma \left(2\,m-2\,j +1\right) \Gamma \left( 2\,j+1 \right) }}\,.
\end{align} 

\item
$\zeta(2m+1)$ expressed in terms of $E_1(\kappa)$ (see \eqref{F8e2}):
\begin{align} \nonumber
\displaystyle { { \left({2}^{2\,m+1} -1 \right) \zeta \left( 2\,m+1 \right) }{}}=&- \left( -1 \right) ^{m}{\pi}^{2\,m}\sum _{j=0}^{m}{\frac {\mathcal{E} \left(2\,j \right) \psi \left(1- 2\,j+2m \right) }{\Gamma \left( 2\,j+1 \right) \Gamma \left( 1-2\,j+2\,m \right) }}
\mbox{}\\ \nonumber
&+{\frac {i}{\pi}\sum _{k=0}^{\infty }{\frac {{E_{1}} \left( i\pi\, \left( k+1/2 \right)  \right) -{E_{1}} \left(-i\pi\, \left( k+1/2 \right)  \right) }{ \left( k
\mbox{}+1/2 \right) ^{2\,m+1}}}}\,.
\end{align}

\item

Bernoulli, Euler and Harmonic numbers are related (see \eqref{Id2a}):

\begin{align} \nonumber
\displaystyle \sum _{j=0}^{m-1}{\frac {\mathcal{E} \left( 2\,j \right) H_{2m-2j} }{\Gamma \left( 1-2\,j+2\,m \right) 
\mbox{}\Gamma \left( 2j+1 \right) }}&=\sum _{k=1}^{m}{\frac {{2}^{2\,k} \left( {2}^{2\,k}-1 \right) B_{2k}  }{ \left( 1-2\,k+2\,m \right) \Gamma \left( 2-2\,k+2\,m \right) 
\mbox{}\Gamma \left( 2k+1 \right) }}\\ \nonumber
&-{\frac {1}{2\,m\Gamma \left( 1+2\,m \right) }}
\end{align}

\item

For integer $p>m,\; m<c_m<m+1$, a contour integral representation for $\zeta(2m+1)$ (see \eqref{Cpa4}):

\begin{align} \nonumber
\displaystyle  &\left( {2}^{2\,m+1}-1 \right) \zeta \left( 2m+1 \right) ={\frac {i \left( -1 \right) ^{p}}{\pi^{2p+1}}\int_{c_m-i\infty }^{c_m+i\infty }\!{\frac {\Gamma \left( 2p-2t \right) {\pi}^{2t}
}{\sin \left( \pi\,t \right) }}\sum _{k=0}^{\infty } \left( -1 \right) ^{k} \left( k+1/2 \right) ^{2\,t-2\,p-2\,m-1}\,{\rm d}t}\\ \nonumber
&- \left( -1 \right) ^{m}{\pi}^{2m}\left(\sum _{j=0}^{p-2-m}{\frac {\Gamma \left( 2\,j+2 \right) \mathcal{E} \left( 2\,j+2+2\,m \right) }{\Gamma \left( 2\,j+3+2\,m \right) }}
\mbox{}+\sum _{j=0}^{m}{\frac {\mathcal{E} \left( 2\,j \right) \psi \left( 1-2\,j+2\,m \right) }{\Gamma \left( 1-2\,j+2\,m \right) \Gamma \left( 2\,j+1 \right) }}\right)\,.
\end{align}

\item

A divergent, asymptotic series representation for the alternating Hurwitz Zeta Function $\eta(1,(t+1)/2)$ (see \eqref{C1Asy}): 
\begin{equation}
\displaystyle \eta(1,(t+1)/2)\sim \sum _{k=0}^{2\,m}{\frac {\mathcal{E} \left( 2\,k \right) }{{t}^{2\,k+1}}}\hspace{20pt}(t\rightarrow\infty), \forall m\geq 0.
\label{AsyInv}
\end{equation}

\item

A representation for the lower incomplete Zeta function $\zeta^{-}(s)$ (see \eqref{ZmEq1}):

\begin{equation} \nonumber
\displaystyle \zeta^{{-}} \left( s \right)={2}^{s-1}\int_{-\pi/2}^{\pi/2}\!{\frac { {{\rm e}^{-\pi\,\sin \left( \theta \right)/2 }}\sin \left( \pi\,\cos \left( \theta \right)/2 +\theta\, \left( s-1 \right)  \right)
\mbox{}}{\cosh \left( \pi\,\sin \left( \theta \right)  \right) -\cos \left( \pi\,\cos \left( \theta \right)  \right) }}\,{\rm d}\theta\,.
\end{equation}

\item

An interesting integral with $0\leq 2n$ (see \eqref{Sa2cEp}):

\begin{equation} \nonumber
\displaystyle \int_{-\pi/2}^{\pi/2}\!{\frac {\sin \left( (\pi/2)\,\cos \left( \theta \right) 
+(2\,n-1)\,\theta \right) {{\rm e}^{-(\pi/2)\,\sin \left( \theta \right) }}
\mbox{}}{\cosh \left( \pi\,\sin \left( \theta \right)  \right) -\cos \left( \pi\,\cos \left( \theta \right)  \right) }}\,{\rm d}\theta=\left( -1 \right) ^{n}{B_{2n\, }{\pi}^{2\,n}\frac {  
\mbox{} \left( {2}^{1-2\,n}-1 \right) }{\Gamma \left( 2\,n+1 \right) }}
\end{equation}
\item

An interesting integral with $2n>0$ (see \eqref{Sa2cEm}):

\begin{equation} \nonumber
\displaystyle \int_{-\pi/2}^{\pi/2}\!{\frac {\sin \left( (\pi/2)\,\cos \left( \theta \right) 
-(2\,n+1)\,\theta \right) {{\rm e}^{-(\pi/2)\,\sin \left( \theta \right) }}
\mbox{}}{\cosh \left( \pi\,\sin \left( \theta \right)  \right) -\cos \left( \pi\,\cos \left( \theta \right)  \right) }}\,{\rm d}\theta=0
\end{equation}

\item
An integral interesting on its own merits (see \eqref{V02p}):

\begin{equation}
\displaystyle \int_{-\pi/2}^{\pi/2}\!{\frac {\cos \left( \theta \right) \sin \left( (\pi/2)\,\cos \left( \theta \right)  \right) {{\rm e}^{-(\pi/2)\,\sin \left( \theta \right) }}}{\cosh \left( \pi\,\sin \left( \theta \right)  \right) 
\mbox{}-\cos \left( \pi\,\cos \left( \theta \right)  \right) }}\,{\rm d}\theta=1/2+{\pi}^{2}/48\,,
\end{equation}
\item

A theorem for $p>0$, (see \eqref{Gans4}):

\begin{align*}
\displaystyle \left( -1 \right) ^{p}\Gamma \left( 2\,p \right)& \frac {    { E_{2p}} \left( -i\pi\, \left( k+1/2 \right)  \right) + {E_{2p}} \left( i\pi\, \left( k+1/2 \right)  \right)  
\mbox{}} { \left( \pi\, \left( k+1/2 \right)  \right) ^{2\,p-1}}\\ 
 &=2\, \left( -1 \right) ^{k}\sum _{j=1}^{p-1}{\frac { \left( -1 \right) ^{j}\Gamma \left( 2\,j \right) }{ \left( \pi\, \left( k+1/2 \right)  \right) ^{2\,j}}} 
- 2 \,(\rm{Si} (\pi  (k + 1/2))-\pi/2 )
\end{align*}

\item

A hypergeometric identity (see \eqref{Corollary1}):

\begin{align*}
 {\mbox{$_1$F$_2$}(\frac{1}{2}-p\,;\,\frac{1}{2},\frac{3}{2}-p\,;\,-\frac{{\pi}^{2}\, \left( k+1/2 \right) ^{2}}{4})}&={\frac { \left( -1 \right) ^{p+1} \left( \pi\, \left( k+1/2 \right)  \right) ^{2\,p}
\mbox{}}{\Gamma \left( 2\,p-1 \right) }}{\mbox{$_1$F$_2$}(\frac{1}{2};\,\frac{3}{2},\frac{3}{2};\frac{-\pi^{2}\,\left( k+1/2 \right)^{2}}{4} )}\\
&+{\frac { \left( -1 \right) ^{k+p} \left( \pi\, \left( k+1/2 \right)  \right) ^{2\,p-1}}{\Gamma \left( 2\,p
\mbox{}-1 \right) }\sum _{j=1}^{p-1}{\frac {\Gamma \left( 2\,j \right)  \left( -1 \right) ^{j}}{ \left( \pi\, \left( k+1/2 \right)  \right) ^{2\,j}}}}
\end{align*}

\end{itemize} 

\section{A series representation for $\zeta(s)$} \label{sec:Series}
By way of review, the standard definition for Dirichlet's function $\eta(s)$ (also called  the alternating zeta function) is given by
\begin{equation}
\eta(s)=\sum_{k=1}^{\infty} (-1)^kk^{-s}
\label{etadef1}
\end{equation}
valid for all $s$, related to $\zeta(s)$ by
\begin{equation}
\eta(s)=(1-2^{1-s})\zeta(s).
\label{etadef2}
\end{equation}
Another well-known representation is
\begin{equation}
(2^s-1)\zeta(s)=\sum_{k=0}^{\infty}1/(k+1/2)^s\,,\hspace{20 pt} \Re(s)>1\,.
\label{ZsumHalf}
\end{equation}
The Cauchy Integral Theorem for $\eta(s)$, with an additional factor $2^{1-s}$ embedded for convenience, is


\begin{equation}
\displaystyle {2}^{1-s} \left( 1-{2}^{1-s} \right) \zeta \left( s \right) ={\frac {1}{2\pi i}\oint \!{\frac { \left( 1-{2}^{1-v} \right) \zeta \left( v \right) {2}^{1-v}
\mbox{}}{s-v}}\,{\rm d}v}
\label{Cauchy}
\end{equation}

where the contour of integration encloses the non-negative real axis as well as the pole at the point $v=s$ in a clockwise direction. Notice that there is no singularity at $v=1$ because of a cancellation between the first two numerator terms of the integrand at $v=1$.  Let $v:=v+1$ in \eqref{Cauchy} giving

\begin{equation}
\displaystyle \eta(s)= \left( 1-{2}^{1-s} \right) \zeta \left( s \right) ={\frac {{2}^{s-1}}{2\pi i}\oint \!{\frac { \left( 1-{2}^{-v} \right) \zeta \left( v+1 \right) {2}^{-v}
\mbox{}}{s-v-1}}\,{\rm d}v}
\label{Cauchy+1}
\end{equation}

where the contour now crosses the real axis at $v<\min(-1,\Re(s))$. Apply the slightly rewritten form of the well-known functional equation
\begin{equation}
\displaystyle \zeta \left( v+1 \right) =\zeta \left( -v \right) {2}^{v+1}{\pi}^{v}\cos \left(\pi\,v/2 \right) \Gamma \left( -v \right) 
\label{Feq5}
\end{equation}
to obtain, with the same contour,

\begin{equation}
\displaystyle  \eta(s)=\left( 1-{2}^{1-s} \right) \zeta \left( s \right) ={\frac {{2}^{s-1}}{i\pi}\oint \!{\frac { \left( 1-{2}^{-v} \right) \zeta \left( -v \right) \Gamma \left( -v \right) {\pi}^{v}
\cos \left( \pi\,v/2 \right) }{s-v-1}}\,{\rm d}v}\,.
\label{Cauchy1}
\end{equation}

Notice that the product $\zeta(-v)\,\Gamma(-v)$ is non-singular on the non-negative real $v-$axis. Substitute the cosine term in \eqref{Cauchy1} as a sum of exponentials giving

\begin{equation}
\eta(s)=\displaystyle  \left( 1-{2}^{1-s} \right) \zeta \left( s \right) =\frac {{2}^{s-1}}{2\pi i}\oint \!{\frac { \left( 1-{2}^{-v} \right) \zeta \left( -v \right)  \left( i\pi \right) ^{v}
\mbox{}\Gamma \left( -v \right) }{s-v-1}}\,{\rm d}v +\{(i\pi)^{v}\rightarrow {(-i\pi)^{ v}} \}.
\label{L3}
\end{equation}

where each exponential has been written as
\begin{equation}
{\rm e}^{\pm i\pi v/2}=(\pm {i})^v\,.
\end{equation}
Now apply \eqref{ZsumHalf} and interchange the summation and integration because for $\Re(v)<-1$ both are convergent, to yield

\begin{equation}
\displaystyle {2}^{1-s} \left( 1-{2}^{1-s} \right) \zeta \left( s \right) =-{\frac {1}{2\pi i}\sum _{k=0}^{\infty } \left( \oint \!{\frac { \left( i\pi \left( k+\frac{1}{2} \right) \right) ^{v}\Gamma \left( -v \right) 
\mbox{}}{s-v-1}}\,{\rm d}v \right) }-{\frac {1}{2\pi i}\sum _{k=0}^{\infty } \left( \oint \!{\frac { \left( -i \pi \left( k+1/2 \right)  \right) ^{v}
\mbox{}\Gamma \left( -v \right) }{s-v-1}}\,{\rm d}v \right) }
\label{L4}
\end{equation}

Because the contour lies to the left of the point $v=s-1$, each of the contour integrals in \eqref{L4} is recognizable as a contour integral representation of the Generalized Integro-Exponential function \cite[Eq. (2.6a)]{Milgram:1985}
\begin{equation}
E_s(z)=\frac{1}{2\pi i}\oint\frac{\Gamma(-v)z^v}{(s-1-v)} {\rm d}v
\label{EiCint}
\end{equation}
where the integration contour encloses the pole at $v=s-1$ and the non-negative $v-$ axis in a clockwise direction. Another alternate (defining) representation is 
\begin{equation}
\displaystyle {\it E_{s}} \left(i\pi\,\kappa \right) \equiv \int_{1}^{\infty }\!{v}^{-s}{{\rm e}^{-i\pi \kappa v}}\,{\rm d}v\,.
\label{ExpInt}
\end{equation}

When \eqref{EiCint} is substituted into \eqref{L4} we find the (slowly converging - see Section \ref{sec:Convergence}) series representation
\begin{equation}
\displaystyle \zeta \left( s \right) ={\frac {{2}^{s-1}
\mbox{}}{{2}^{1-s}-1}}\sum _{k=0}^{\infty }\left({\it E_{s}} \left(i\pi\, \left( k+1/2 \right)  \right) +{\it E_{s}} \left(-i\pi\, \left( k+1/2 \right)  \right)\right) \,. 
\label{EiSum}
\end{equation}

{\bf Remark:} Because it has never been proven (\cite[Section 2.13]{Titch2}) that the functional equation \eqref{Feq5} {\it uniquely} defines $\zeta(s)$ (and therefore $\eta(s)$), up to the point in \eqref{L4} where \eqref{ZsumHalf} is applied, \eqref{L3} is valid for any function (that could have been labelled $\zeta(s)$ arbitrarily) which might obey \eqref{Feq5} in any of its forms.

\subsection{Convergence}\label{sec:Convergence}


Although \eqref{EiSum} converges slowly numerically, it does converge. To test this observation, define

\begin{equation}
\displaystyle t^{\pm}_{{k}} \left( s \right) ={ E_{s}} \left(i\pi\,(k+1/2) \right) \pm { E_{s}} \left(- i\pi\,(k+1/2) \right) 
\label{tkp(s}
\end{equation}

As $k\rightarrow \infty$
\begin{equation}
\displaystyle {\it t^{+}_{k}} \left( s \right)\sim  \left( -1 \right) ^{k}\left( -{\frac {2}{\pi\,k}}+{\frac {1}{\pi\,{k}^{2}}}+{\frac {1}{\pi\,{k}^{3}} \left( -\frac{1}{2}+{\frac {2s \left( s+1 \right) }{{\pi}^{2}}} \right) } \right) +O \left( {k}^{-4} \right) \,.
\label{Tk(s)Inf}
\end{equation} 
Applying Gauss' test gives

\begin{equation}
\displaystyle  \left| {\frac {t^{+}_{{k+1}} \left( s \right) }{t^{+}_{{k}} \left( s \right) }} \right| \sim 1-\frac{1}{k}+\frac{3}{2k^2}+O \left( {k}^{-3} \right) 
\label{tkpRat}
\end{equation}
indicating ambiguity. However, from \eqref{Tk(s)Inf}, we see that, for finite values of $s$, $|t^{+}_{k}(s)|\rightarrow 0$ monotonically, and therefore, by the Alternating Series test, the sum converges. 

In the case of ${\it t^{-}_{k}}(s)$ we have

\begin{equation}
\displaystyle {\it t^{-}_{k}} \left( s \right)\sim{\frac {is \left( -1 \right) ^{k}}{{\pi}^{2}} \left(-\frac{2}{k^2} +\frac{2}{k^3}+{\frac {1}{{k}^{4}} \left(-\frac{3}{2} +{\frac { 2\,\left( s+1 \right)  \left( s+2 \right) }{{\pi}^{2}}} \right) }
\mbox{} \right) }
\label{tkm(s}
\end{equation}
and, with respect to Gauss' test,
\begin{equation}
\displaystyle  \left| {\frac {{\it t^{-}}_{{k+1}} \left( s \right) }{{\it t^{-}}_{{k}} \left( s \right) }} \right| =1-\frac{2}{k}+O \left( {k}^{-2} \right) 
\label{tkmRat}
\end{equation} 
which proves convergence.\newline

In a later Section, we will consider sums involving $t^{\pm}_{k}(s+p)/(k+1/2)^{p}$ where $p>0$. Testing such sums for $p:=2p$, gives

\begin{equation}
\displaystyle  \left| {\frac {t^{+}_{{k+1}} \left( s+2\,p \right) }{t^{+}_{{k}} \left( s+2\,p \right) }} \right| \sim 1-\frac{1+2p}{k}\,,
\label{tkRatEven}
\end{equation}
proving absolute convergence by Gauss' test, and for $p:=2p-1$ 

\begin{equation}
\displaystyle  \left| {\frac {t^{-}_{{k+1}} \left( s+2p-1 \right) }{t^{-}_{{k}} \left( s+2p-1 \right) }} \right| \sim 1-\frac{2p}{k}\,,
\label{tkRatOdd}
\end{equation}
again proving absolute convergence by the same test.

\subsection{Recursion} \label{sec:recursion}
The recursion rule for $E_s(z)$ is
\begin{equation}
\displaystyle {\it E_s} \left(z \right) ={\frac {z E_{s-1} \left(z \right) -{{\rm e}^{-z}}}{1-s}}\hspace{15pt}s\neq1,
\label{recur}
\end{equation}

leading, after the use of forward recursion, to  

\begin{equation}
\displaystyle \eta(s)=(1-{2}^{1-s})\zeta \left( s \right) ={2}^{s-1}\,\frac{is}{\pi}{\frac {}{ }\sum _{k=0}^{\infty }{\frac {{\it E}_{s+1} \left(-i\pi\, \left( k+1/2 \right)  \right) -{\it E}_{s+1} \left(i\pi\, \left( k+1/2 \right)  \right) 
\mbox{}}{k+1/2}}}+{ {{2}^{s-1}}}
\label{R2a}
\end{equation}
with improved convergence properties (see \eqref{tkRatEven}). The general form of \eqref{R2a} can be obtained similarly by induction, giving for any $n\geq 0$


\begin{align} \nonumber
\displaystyle \eta(s)=&{ \left(1-{2}^{1-s} \right) \zeta \left( s \right)  }= \frac{{2}^{s}}{\Gamma \left( s \right)}\sum _{j=0}^{\left\lfloor (n-1)/2\right\rfloor}  \frac{\left( -1 \right) ^{j}
\mbox{}\Gamma \left( s+2\,j \right)}{\pi^{2j+1}} \sum _{k=0}^{\infty }{\frac { \left( -1 \right) ^{k}}{  \left( k+1/2 \right)  ^{2\,j+1}}} \\
&-\frac{{2}^{s-1}\,\Gamma \left( n+s \right)}{\Gamma \left( s \right)}\left( \frac{-1}{\;\pi} \right) ^{n}{{\rm e}^{i\pi\,n/2}}\sum _{k=0}^{\infty }\left({\frac {{E_{n+s}} \left( -i\pi\, \left( k+1/2 \right)  \right)+ \left( -1 \right) ^{n}{E_{n+s}} \left( i\pi\, \left( k+1/2 \right)  \right)  }{ \left( k+1/2 \right) ^{n}}}\right)\,.
\label{R5}
\end{align}
As $n$ increases in \eqref{R5}, the two terms composing the right-hand side converge to each other, so numerically significant error occurs at large values of $n$, although the convergence rate of the second infinite sum improves also. An equivalent form may be obtained by replacing $s:=1-s$ and applying the functional equation \eqref{Feq5} to obtain

\begin{align} \nonumber
\displaystyle & \left( {2}^{s}-1 \right) \zeta \left( s \right) =-{\pi}^{s-1}\sin \left( \pi s/2 \right)\sum _{j=0}^{\left\lfloor (n-1)/2\right\rfloor}{\frac {\Gamma \left( 2\,j+1-s \right) 
\mbox{}{\mathcal{E}(2j)} }{\Gamma \left( 2\,j+1 \right) }} \\
&+{\frac {{\pi}^{s-n}{{\rm e}^{i\pi\,n/2}}}{2\Gamma \left( s-n \right) \cos \left( \pi s/2 \right) }\sum _{k=0}^{\infty }{\frac {{E_{1+n-s}} \left(-i\pi\, \left( k+1/2 \right)  \right) + \left( -1 \right) ^{n}{E_{1+n-s}} \left(i\pi\, \left( k+1/2 \right)  \right) 
\mbox{}}{ \left( k+1/2 \right) ^{n}}}}
\label{R5b}
\end{align}
In \eqref{R5b}, the infinite inner sum  appearing in the first term of \eqref{R5} has been replaced by a known result in terms of Euler numbers $\mathcal{E}(2j)$ - see \eqref{BetaDef}. The application of the result \eqref{R5} is fairly straightforward for both even and odd values of $n$ if $s\neq \pm m$. These exceptional cases are investigated in the following subsections. 


\subsubsection{Eq. \eqref{R5} with $n$ even or odd yields a recursion formula for Bernoulli numbers} \label{sec:nEven}
Consider the case $n:=2n$. Application of the reversed logic presented in Section \ref{sec:Series} applied to \eqref{R5} along with \eqref{ExpInt} leads to

\begin{equation}
\displaystyle \sum _{k=0}^{\infty }{\frac {{E_{s+2n}} \left(-i\pi\, \left( k+1/2 \right)  \right) +{E_{s+2n}} \left( i\pi\, \left( k+1/2 \right)  \right) 
\mbox{}}{ \left( k+1/2 \right) ^{2\,n}}}=2\,\int_{1}^{\infty }\!\sum _{k=0}^{\infty }{\frac {{v}^{-s-2\,n}\cos \left( \pi\, \left( k+1/2 \right) v \right) }{ \left( k
\mbox{}+1/2 \right) ^{2\,n}}}\,{\rm d}v\,.
\label{Es1a}
\end{equation}

Expand the integration
\begin{equation} 
\int_{1}^{\infty}...{\rm d}v=\int_{0}^{\infty}...{\rm d}v-\int_{0}^{1}...{\rm d}v
\label{IntSplit}
\end{equation}

interchange the (convergent) sum and integral, and consider the first integration term on the right of \eqref{IntSplit}


\begin{align} \nonumber
\displaystyle \sum _{k=0}^{\infty }  \int_{0}^{\infty }\!{\frac {{v}^{-s-2\,n}\cos \left( \pi\, \left( k+1/2 \right) v \right) }{ \left( k+1/2
\mbox{} \right) ^{2\,n}}}& \,{\rm d}v  =\left( -1 \right) ^{n}{\pi}^{s+2\,n-1}\sum _{k=0}^{\infty } \left( k+1/2 \right) ^{s-1}\Gamma \left(1 -s-2\,n \right) \sin \left( \pi s/2 \right) \\ \nonumber
 = &\left( -1 \right) ^{n}{\pi}^{s+2n-1} \left( {2}^{1-s}-1 \right) \zeta \left( 1-s \right) \Gamma \left( 1-s-2n \right) 
\sin \left( \pi s/2 \right)\\ 
=&\displaystyle {\frac {{\pi}^{2\,n} \left({2}^{1-s} -1 \right) \zeta \left( s \right)  \left( -1 \right) ^{n}\Gamma \left( s \right) 
\mbox{}}{{2}^{s}\Gamma \left( s+2\,n \right) }} \,.
\label{FullInt}
\end{align}

The first equality in \eqref{FullInt} follows from the application of \eqref{Erd37} with $-1<-\Re(s)-2n<0$, the second equality follows from the first by the principle of analytic continuation using \eqref{ZsumHalf} with $\Re(s)<0$, and the third equality follows from simple identities and invocation of the functional equation \eqref{Feq5} with $v:=v-1$.\newline

With regards to the second integral term in \eqref{IntSplit}, apply \eqref{hans1} without interchanging the sum and integration, to find

\begin{align} \nonumber
\displaystyle \int_{0}^{1}\!{v}^{-s-2\,n}&\sum _{k=0}^{\infty }{\frac {\cos \left( \pi\, \left( k+1/2 \right) v \right) }{ \left( k+1/2 \right) ^{2\,n}}}\,{\rm d}v=
\mbox{}\,{\frac { \left( 2\,\pi \right) ^{2\,n} \left( -1 \right) ^{n}}{4\,\Gamma
\mbox{} \left( 2\,n \right) }\int_{0}^{1}\!{v}^{-s-2\,n}{\mathcal{E}_{ 2\,n-1}} \left(v/2 \right) \,{\rm d}v}\\
&= \left( -1 \right) ^{n+1}\,\frac{\left( 2\,\pi \right) ^{2\,n} }{2}\sum _{k=0}^{2\,n-1}{\frac {  \left( {2}^{-k}-{2}^{-2\,k+2\,n} \right)B_{2n-k} }{\Gamma \left( k+1 \right) \Gamma \left(1 -k+2\,n \right) 
\mbox{} \left( -1+s+2\,n-k \right) }}
\label{Int2}
\end{align}

The second equality in \eqref{Int2} arises from the definition \eqref{EulernoRev} followed by straightforward integration with $-s-2n>-1$. Now apply \eqref{Es1a},  \eqref{IntSplit}, \eqref{FullInt} and \eqref{Int2} into \eqref{R5}, along with the use of  \eqref{BetaDef} to identify the last (inner) sum of \eqref{R5} in terms of Euler numbers, to arrive at the identity


\begin{equation}
\displaystyle \sum _{j=0}^{n-1}{\frac {\Gamma \left( 2\,j+s \right) {\mathcal{E}} \left( 2\,j \right) }{\Gamma \left( 2\,j+1 \right) }}= {2}^{2\,n}\Gamma \left( s+2\,n \right) \,\sum _{k=0}^{2\,n-1}\frac {  \left( {2}^{-k}-{2}^{2\,n-2\,k} \right){B_{2n-k}} }{\Gamma \left( k+1 \right) \Gamma \left(1 -k+2\,n \right) 
\mbox{} \left( 2\,n-1+s-k \right) }\,.
\label{EuId}
\end{equation}

This result can be used to obtain an interesting recursion relation among Bernoulli numbers $B_{n}$ due to the presence of the arbitrary parameter $s$. To solve \eqref{EuId}, first let $n:=N$, then let $n:=N+1$ and subtract, giving


\begin{align} \nonumber
\displaystyle \frac {\mathcal{E} \left( 2\,N \right) }{\Gamma \left( 2\,N+1 \right) }&=-{2}^{2\,N}\sum _{k=0}^{2\,N-1}{\frac {{B_{2N-k}}   \left( {2}^{-k}-{2}^{-2\,k+2\,N} \right) }{\Gamma \left( k+1 \right) \Gamma \left( -k+2\,N+1 \right) 
\mbox{} \left( -1+s+2\,N-k \right) }}\\&+{\frac {\Gamma \left( s+2\,N+2 \right) {2}^{2\,N+2}}{\Gamma \left( s+2\,N \right) }\sum _{k=0}^{2\,N+1}{\frac {{B_{2-k+2N}}  \left( {2}^{-k}-{2}^{2-2\,k+2\,N} \right) }{\Gamma \left( k+1 \right) \Gamma \left(3 -k+2\,N \right) 
 \left( 1+s+2\,N-k \right) }}}\,,
\label{EuId2}
\end{align}

which, after shifting the index of the second sum by two, reorganizing, and recognizing that $B_{2N+1}=0$ (except define $B_{1}=-1/2$), yields a result of the following form:

\begin{equation}
\displaystyle \frac {\mathcal{E} \left( 2\,N \right) }{\Gamma \left( 2\,N+1 \right) }=X(N)s+Y(N)\,,
\label{Sdep}
\end{equation}
where $X(N)$ and $Y(N)$ are complicated functions only of the variable $N$. Because the left-hand side of \eqref{Sdep} does not depend on $s$, it must be that $X(N)=0$, and when $X(N)=0$ is simplified, one obtains

\begin{equation}
\displaystyle {B_{2N+2}}  ={\frac {\Gamma \left( 3+2\,N \right) }{{2}^{2\,N} \left( {2}^{2\,N+2}-1 \right) }\sum _{k=0}^{N-1}{\frac {{2}^{2\,k}  \left( 1-{2}^{2\,k+2} \right)B_{2k+2} }{\Gamma \left( 1-2\,k+2\,N \right) \Gamma \left( 2\,k+3 \right) 
\mbox{}}}}+{\frac {1+N}{{2}^{2\,N+1} \left( {2}^{2\,N+2}-1 \right) }}, \hspace{20pt} N\geq0\,, 
\label{BernId}
\end{equation} 

one of a large number of possible recursive definitions of Bernoulli numbers $B_{N}$ (e.g. \cite[Eqs. (4) and (6)]{Tyagi}, \cite{Omer} and compare with \cite[Eq. 24.5.8]{NIST}). This result can alternatively be derived from \eqref{EuNo} using the fact that $\mathcal{E}(2n+1)=0$ as well as from a convolution type of identity for Bernoulli polynomials at $x=1/2$ and $x=1/4$, the latter (privately) communicated by C. Vignat, Tulane University. The remaining term $Y(N)$ of \eqref{Sdep} eventually simplifies to reproduce \eqref{EuNo}.
\newline

By following steps similar to those given above, the same results emerge for the case $n:=2n+1$.

\subsection{Eq. \eqref{R5b} with $s=2m$} \label{sec:s_2m}

The evaluation of \eqref{R5b} with either $s=2m$ or $s=2m+1$ requires that limits be taken. Consider the case $s=2m$. Substitute $s=2m+\epsilon$ in \eqref{R5b} and evaluate the limit $\epsilon\rightarrow0$ to obtain  


\begin{equation}
\displaystyle  \left( 1-{2}^{2\,m} \right) \zeta \left( 2\,m \right) = \left( -1 \right) ^{m}\frac{{\pi}^{2\,m}}{2}\sum _{j=0}^{m-1}{\frac {\mathcal{E} \left( 2\,j \right) }{\Gamma \left( 2\,m-2\,j \right) \Gamma \left( 2\,j+1 \right) }}\,,
\label{CgId2}
\end{equation}
a result equivalent to \cite[Eq. 2.4.15]{NIST}, valid for all non-negative values of the arbitrary parameter $n$. Also, note both \cite[Eq. (3.1)]{milgram2020additions} as well as \cite[Eq. 24.5.4]{NIST} where we find the related identity

\begin{equation}
\displaystyle \sum _{j=0}^{m}{\frac {\mathcal{E} \left( 2\,j \right) }{ \Gamma \left( 2\,m-2\,j+1 \right)
\Gamma \left( 2\,j+1 \right) }}=0\,.
\label{EuSum}
\end{equation}

\subsection{Eq. \eqref{R5b} with $s=2m+1$, $n<2m+1$, and $m\geq0$} \label{sec:s_2mP1_Lt}


For the case $s=2m+1$, again set $s=2m+1+\epsilon$ and evaluate the limit $\epsilon\rightarrow0$ with $n<2m+1$. In this case, on the right-hand side, a term of order $\epsilon^{-1}$ arises, the coefficient of which must vanish (because the left-hand side is not divergent) giving rise to an identity as follows. Let $n=2m-p$ with $0\leq p\leq 2m$ and use \eqref{EuSum} to find

\begin{align} \nonumber
\displaystyle \frac{{\pi}^{p-2\,m} \left( -1 \right) ^{m}{{\rm e}^{-i\pi p/2}}}{\Gamma \left( 1+p \right)}
\mbox{}\sum _{k=0}^{\infty }{\frac {{ \left( -1 \right) ^{p}\,E_{-p}} \left( i\pi\, \left( k+1/2 \right)  \right) 
\mbox{}+{E_{-p}} \left(-i\pi\, \left( k+1/2 \right)  \right) }{ \left(k+1/2 \right) ^{-p+2\,m}}}\\=- \sum _{j=\left\lfloor -p/2+m+1/2\right\rfloor}^{m}{\frac {\mathcal{E}
 \left( 2\,j \right) }{\Gamma \left(2\,m-2\,j +1\right) \Gamma \left( 2\,j+1 \right) }}\,.
\label{T0aA}
\end{align} 


Notice that if $p=2r$ or $p=2r+1$, the right-hand side of \eqref{T0aA} is invariant for constant values of $m$ and therefore so is the left-hand side, leading to some interesting identities. Suppose $p=2m$, in which case, the right-hand side of \eqref{T0aA} vanishes by \eqref{EuSum}, giving

\begin{equation}
\displaystyle \sum _{k=0}^{\infty }{E_{-2m}} \left( i\pi\, \left( k+1/2 \right)  \right) +{E_{-2m}} \left( -i\pi\, \left( k+1/2 \right)  \right) 
\mbox{}=0\,.
\label{L2m}
\end{equation}
Further, if $p=2m-1$, by extending the lower limit of the finite sum in \eqref{T0aA}, we similarly find
\begin{equation}
\displaystyle \sum _{k=0}^{\infty }{\frac {{E_{1-2m}} \left( i\pi\, \left( k+1/2 \right)  \right) -{E_{1-2m}} \left(-i\pi\, \left(k+1/2 \right)  \right) 
\mbox{}}{k+1/2}}={\frac {i\pi}{2m}}\,,
\label{L3a}
\end{equation}
and, if $p=2m-2$,
\begin{equation}
\displaystyle \sum _{k=0}^{\infty }{\frac {{E_{2-2m}} \left(i\pi\, \left( k+1/2 \right)  \right) +{E_{2-2m}} \left(-i\pi\, \left( k+1/2 \right)  \right) 
\mbox{}}{ \left( k+1/2 \right) ^{2}}}=-{\frac {{\pi}^{2}\Gamma \left( 2\,m-1 \right) }{\Gamma \left( 2\,m+1 \right) }}
\label{L4a}
\end{equation}

{\bf Remark.} Replacing ${E_{-p}} \left(\pm\,  i\pi\, \left( k+1/2 \right)  \right)$ in \eqref{T0aA} by the finite sum \eqref{EiNegP}, reduces it to a tautology; this replacement also provides an alternative proof of \eqref{L2m}, \eqref{L3a} and \eqref{L4a} and their obvious extensions. Furthermore, setting $s=-2m$ in \eqref{EiSum} with recourse to \eqref{L2m} yields the well-know result $\zeta(-2m)=0$ - the trivial zeros. Finally, although \eqref{L3a} and \eqref{L4a}, related by $p=2r$ and $p=2r+1$ with $r=m-1$ appear to differ (in violation of the remark immediately following \eqref{T0aA}), that is only because $p$-dependent factors have been transposed to the right-hand side for the sake of uniformity of presentation. \newline

Consider now the terms associated with $\epsilon^{0}$. The analysis is long and complicated, involving the derivative of $E_s(z)$ with respect to its order $s$ (see \cite{Milgram:1985}) as well as considerable simplification, all of which employ identities listed in \cite[Eqs. (2.1), (2.20c), and (2.21)]{Milgram:1985}. As before, the result is independent of the parameter $n$, so, a simple instance of the final result with $n=0$ is


\begin{align} \nonumber
\displaystyle { { \left({2}^{2\,m+1} -1 \right) \zeta \left( 2\,m+1 \right) }{}}=&- \left( -1 \right) ^{m}{\pi}^{2\,m}\sum _{j=0}^{m}{\frac {\mathcal{E} \left(2\,j \right) \psi \left(1- 2\,j+2m \right) }{\Gamma \left( 2\,j+1 \right) \Gamma \left( 1-2\,j+2\,m \right) }}
\mbox{}\\ 
&+{\frac {i}{\pi}\sum _{k=0}^{\infty }{\frac {{E_{1}} \left( i\pi\, \left( k+1/2 \right)  \right) -{E_{1}} \left(-i\pi\, \left( k+1/2 \right)  \right) }{ \left( k
\mbox{}+1/2 \right) ^{2\,m+1}}}}\,.
\label{F8e2}
\end{align}
For other values of $n$, the equivalent of \eqref{F8e2} reduces to a tautology. Notice that both sides of \eqref{F8e2} are singular if $m=0$ because of \eqref{L3a}.

\subsection{Eq. \eqref{R5b} with $s=2m+1$ and $n> 2m$} \label{sec:s_2mP1Gt}

For the same case ($s=2m+1$) this time with $n>2m$, in order to repeat the same procedure as above, it is necessary to split the finite sum in \eqref{R5b} into two parts: $0\leq j\leq m$ and $m+1\leq j\leq \lfloor(n-1)/2\rfloor $. The numerator and denominator of the first of these sums both vanish when $s\rightarrow2m+1$ so a limit must be evaluated, giving

\begin{equation}
\displaystyle \lim _{s\rightarrow 2\,m+1} \left( {
 \frac{\displaystyle \sum _{j=0}^{m}{\frac {\mathcal{E} \left( 2\,j \right) }{\Gamma \left( s-2\,j \right) \Gamma \left( 2\,j+1 \right) }}} {\cos
 \left( \pi\,s/2 \right)}} \right)
  = \left( -1
 \right) ^{m}{\frac {2}{\pi}\sum _{j=0}^{m}{\frac {\mathcal{E} \left( 2\,j \right) \psi \left( 1-2\,j+2\,m \right) }{\Gamma \left( 1-2\,j+2\,m \right) \Gamma \left( 2\,j+1 \right) }}}\,.
 \label{Tlimit}
\end{equation}
Set $n=2m+p,\, p>0$, and after some simplification \eqref{R5b} becomes


\begin{align} \nonumber
\displaystyle  \left( {2}^{2\,m+1}-1 \right) \zeta \left( 2\,m+1 \right) = \frac{\left( -1 \right) ^{p}\Gamma \left( p \right) {{\rm e}^{i\pi p/2}}}{\pi^{p}}\sum _{k=0}^{\infty }{\frac {{E_p} \left(-i\pi\, \left( k+1/2 \right)  \right) + \left( -1 \right) ^{p}{E_p} \left( i\pi\, \left( k+1/2 \right)  \right) 
}{ \left( k+1/2 \right) ^{p+2\,m}}}
\mbox{} \\ 
- \left( -1 \right) ^{m}{\pi}^{2\,m}\sum _{j=m+1}^{\left\lfloor p/2+m+1/2\right\rfloor-1}{\frac {\Gamma \left( 2\,j-2\,m \right) 
\mbox{}\mathcal{E} \left( 2\,j \right) }{\Gamma \left( 2\,j+1 \right) }}
\mbox{} 
- \left( -1 \right) ^{m}{\pi}^{2\,m}\sum _{j=0}^{m}{\frac {\mathcal{E} \left( 2\,j \right) \psi \left( 1-2\,j+2\,m \right) }{\Gamma \left( 1-2\,j+2\,m \right) \Gamma \left( 2\,j+1 \right) }}
\label{Casep2a}
\end{align}
Notice that the left-hand side of \eqref{Casep2a} is independent of $p$ and therefore so must be the right-hand side. 

\subsection{Euler, Bernoulli and Harmonic numbers are related}

An interesting identity emerges from \eqref{Casep2a} by applying \eqref{ExpInt} to the first term on its right-hand side with $p=1$ and identifying the sine integral ${Si}$ as follows

\begin{equation}
\displaystyle  \int_{1}^{\infty }\!{\frac {\sin \left( \pi\, \left(j+1/2 \right) v \right) }{v}}\,{\rm d}v=-{\it Si} \left( \pi\,(j+1/2) \right) +\pi/2
 =\pi/2-\int_{0}^{1}\!{\frac {\sin \left( \pi\, \left( j+1/2 \right) v \right) }{v}}\,{\rm d}v\,.
\label{sint}
\end{equation}

Each of the terms on the right-hand side of the second equality of \eqref{sint} can be summed over the index $j$ according to \eqref{T0aA} with recourse to \eqref{ZsumHalf}, \eqref{Hans2}, \eqref{EulerNo} and straightforward integration, in that order ({\bf Remark:} the term corresponding to $j=0$ vanishes because, for $m>0$, $B_{2m+1}=0$), to eventually yield the identity

\begin{equation}
\displaystyle \sum _{j=0}^{m}{\frac {\mathcal{E} \left( 2\,j \right) \psi \left( 1-2\,j+2\,m \right) }{\Gamma \left( 1-2\,j+2\,m \right) 
\mbox{}\Gamma \left( 2\,j+1 \right) }}=\sum _{j=0}^{2\,m-1}{\frac { \left( {2}^{2+2\,j}-{2}^{j+1} \right) B_{j+1}  }{ \left( 2\,m-j \right) 
\mbox{}\Gamma \left( 1-j+2\,m \right) \Gamma \left( j+2 \right) }}\,.
\label{EuBernId}
\end{equation}

Also see \eqref{F8e2}. Since $\psi(1)=-\gamma$, \eqref{EuBernId} is intriguing because it suggests the existence of an algebraic relationship between the Euler-Mascheroni constant $\gamma$, whose rationality status is unknown, and the irrational number $\pi$ embedded in the Bernoulli numbers $B_{j+1}$. This is not the case however.\newline

{\bf Proof:}

From its recursion property
\begin{equation}
\displaystyle \psi \left( 1-2\,j+2\,m \right) =\sum _{k=0}^{2\,m-2\,j-1} \frac{1}{\left( 2\,m-2\,j-k \right)}-\gamma\,,
\label{Psip}
\end{equation}

substitute \eqref{Psip} into \eqref{EuBernId}, reverse the inner sum and separate all terms containing $\gamma$ to find

\begin{align} \nonumber
\displaystyle \sum _{j=0}^{m} {\frac {\displaystyle \mathcal{E} \left( 2\,j \right) {\sum _{k=0}^{2\,m-2\,j-1} \frac{1}{ k+1 }}}{\Gamma \left( 1-2\,j+2\,m \right) 
\mbox{}\Gamma \left( 2\,j+1 \right) }} -&\gamma\,\sum _{j=0}^{m}{\frac {\mathcal{E} \left( 2\,j \right) }{\Gamma \left( 1-2\,j+2\,m \right) \Gamma \left( 2\,j+1 \right) }}\\
&=\sum _{k=0}^{2\,m-1}{\frac { \left( {2}^{2+2\,k}-{2}^{k+1} \right) 
\mbox{}{B_{k+1}} }{ \left( 2\,m-k \right) \Gamma \left( 1-k+2\,m \right) \Gamma \left( k+2 \right) }}\,.
\label{GamTerm}
\end{align}

Note that the coefficient of the $\gamma$ term in \eqref{GamTerm} vanishes by \eqref{EuSum} and therefore, effectively, $\gamma$ is not embedded in \eqref{EuBernId}. {\bf QED}\newline

From \eqref{GamTerm}, after setting $k:=2k+1$ on the right-hand side we are left with an interesting identity that relates Euler, Bernoulli and Harmonic numbers:

\begin{equation}
\displaystyle \sum _{j=0}^{m-1}{\frac {\mathcal{E} \left( 2\,j \right) H_{2m-2j} }{\Gamma \left( 1-2\,j+2\,m \right) 
\mbox{}\Gamma \left( 2j+1 \right) }}=\sum _{k=1}^{m}{\frac {{2}^{2\,k} \left( {2}^{2\,k}-1 \right) B_{2k}  }{ \left( 1-2\,k+2\,m \right) \Gamma \left( 2-2\,k+2\,m \right) 
\mbox{}\Gamma \left( 2\,k+1 \right) }}-{\frac {1}{2\,m\Gamma \left( 1+2\,m \right) }}
\label{Id2a}
\end{equation}

This result can be used to obtain a recursive formula for $\mathcal{E}(2j)$ by inserting \eqref{B2n} into \eqref{Id2a} after setting $m:=m+1$, and then considering the final term ($k=m$) in each sum to eventually obtain


\begin{align}
\displaystyle \mathcal{E} \left( 2\,m \right) =  -4\,\sum _{k=0}^{m-1} {\frac {\mathcal{E} \left( 2
\,k \right) }{\Gamma \left( 2\,k+1 \right) }\sum _{j=0}^{m-k}{\frac { \Gamma \left( 1+2\,m \right)}{\Gamma \left( 2\,m+2-2\,k-2\,j \right)  \left( 2\,j+1 \right) \Gamma \left( 2+2\,j \right) }}} \\ 
 +4\, \Gamma \left( 1+2\,m \right)\sum _{k=0}^{m-1}  {\frac {\mathcal{E} \left( 2\,k \right) H_{2m-2k+2}\,}{\Gamma \left( 3-2\,k+2\,m \right) \Gamma \left( 2\,k+1 \right) }}  
 +{\frac {1}{ \left( 1+2\,m \right)  \left( m+1 \right) ^{2}}}\,.
\label{Rc1b}
\end{align}

Further, by substituting the integral representation of $H_{2m-2j}$ (see \cite[first entry, with $t\rightarrow 1-x$]{HarmonicInt}) into the left-hand side of \eqref{Id2a}, inverting the summation and integration operators, making a change of integration variables ($x\rightarrow 2x$) and utilizing \cite[Eq. 24.2.10]{NIST} to evaluate the sum, we arrive at


\begin{equation}
\displaystyle {\frac {{2}^{2\,m}}{\Gamma \left( 2\,m+1 \right) }\int_{0}^{1}\!{\frac {\mathcal{E} \left( 2\,m,u/2 \right) }{u}}\,{\rm d}u}=-\sum _{j=0}^{m-1}{\frac {\mathcal{E} \left( 2\,j \right) {H} _{ 2\,m-2\,j } 
\mbox{}}{\Gamma \left( 2\,j+1 \right) \Gamma \left( 2\,m-2\,j+1 \right) }}\hspace{30pt} m>0,
\label{EuInt}
\end{equation}

an addition to the collection of similar integrals evaluated in \cite[Section 10]{MollVignat}. {\bf Remark:} The right-hand side of \eqref{Id2a} can be also be reduced to the left-hand side of \eqref{EuInt}, thereby providing an alternative derivation of \eqref{Id2a} - (C. Vignat, Toulane University, private communication.)

\subsection{Contour and various integral representations for $\zeta(2m+1)$}

Further interesting results can be found by considering \eqref{Casep2a} with $n=2m+p$ as before, then setting $p:=2p$ giving
\begin{align} \nonumber
\displaystyle \left( {2}^{2\,m+1}-1 \right) \zeta \left( 2\,m+1 \right)& + \left( -1 \right) ^{m}{\pi}^{2\,m}\sum _{j=0}^{m}{\frac {\mathcal{E} \left( 2\,j \right) \psi \left( 1-2\,j+2\,m \right) }{\Gamma \left( 1-2\,j+2\,m \right) 
\mbox{}\Gamma \left( 2\,j+1 \right) }} \\ \nonumber
&=\left( -1 \right) ^{p}\frac{\Gamma \left( 2\,p \right)}{ {\pi}^{2\,p}} \sum _{k=0}^{\infty }{\frac {{E_{2p}} \left( -i\pi\, \left( k+1/2 \right)  \right) +{E_{2p}} \left(i\pi\, \left( k+1/2 \right)  \right) }{ \left( k
\mbox{}+1/2 \right) ^{2\,p+2\,m}}}\\
&- \left( -1 \right) ^{m}{\pi}^{2\,m}\sum _{j=0}^{p-2}{\frac {\Gamma \left( 2\,j+2 \right) \mathcal{E} \left( 2\,j+2+2\,m \right) }{\Gamma \left( 2\,j+3+2\,m \right) }}\,.
\label{Cp}
\end{align}


The result \eqref{Cp} can be analyzed with recourse to Meijer's G-function; since the analysis is somewhat lengthy the details have been relegated to Appendix \eqref{sec:Circuitous}. With reference to that Appendix, a simple change of variables in \eqref{Cpa4}, eventually yields
\begin{align} 
\displaystyle  \left( {2}^{2\,m+1}-1 \right) \zeta \left( 2\,m+1 \right) =A(m)
 - \left( -1 \right) ^{m}{\pi}^{2\,m}\sum _{j=0}^{m}{\frac {\mathcal{E} \left( 2\,j \right) \psi \left( 1-2\,j+2\,m \right) }{\Gamma \left( 1-2\,j+2\,m \right) \Gamma \left( 2\,j+1 \right) 
\mbox{}}}\,.
\label{Cp5a}
\end{align}
where 
\begin{equation}
A(m)\equiv {\frac {1}{{\pi}^{2}}\int_{-\infty }^{\infty }\!{\frac {\Gamma \left( 1-2\,iv \right) 
\mbox{}{\pi}^{2\,iv}\sum _{k=0}^{\infty } \left( -1 \right) ^{k} \left( k+1/2 \right) ^{-2\,m-2+2\,iv}
\mbox{}}{\cosh \left( \pi\,v \right) }}\,{\rm d}v}
\label{Adef}
\end{equation}


The sum in the integrand of $A(m)$ has a known integral representation \cite[Eq. 3.523(3)]{G&R}:
\begin{equation}
\displaystyle \sum _{k=0}^{\infty }  \frac{\left( -1 \right) ^{k} }{ \left( k+1/2 \right)  ^{2\,m+2-2\,iv}}=
\mbox{}\,{\frac {{2}^{2\,m+1-2\,iv
\mbox{}}}{\Gamma \left( 2\,m+2-2\,iv \right) }\int_{0}^{\infty }\!{\frac {{t}^{2\,m+1-2\,iv}}{\cosh \left( t \right) }}\,{\rm d}t}\,.
\label{BetaSum}
\end{equation}
Applying \eqref{BetaSum} to \eqref{Adef} and reordering the integration operators gives

\begin{equation}
A(m)=\displaystyle {\frac {{2}^{2\,m+1}}{{\pi}^{2}}\int_{0}^{\infty }\!{\frac {{t}^{2\,m+1}}{\cosh \left( t \right) }\int_{-\infty }^{\infty }\!{\frac {\Gamma
\mbox{} \left( 1-2\,iv \right) }{\Gamma
\mbox{} \left( 2\,m+2-2\,iv \right) 
\mbox{}\cosh \left( \pi\,v \right) } \left( {\frac {\pi}{2t}} \right) ^{2\,iv}}\,{\rm d}v}\,{\rm d}t}\,.
\label{Cp4}
\end{equation}

Since the arguments of the $\Gamma$ function ratio in \eqref{Cp4} differ by an integer, this ratio can profitably be written as a $\beta$ function (not Dirichlet's Beta function), which in turn has a well-known integral representation, giving
\begin{equation}
\displaystyle {\frac {\Gamma \left( 1-2\,iv \right) }{\Gamma \left( 2\,m+2-2\,iv \right) }}=  {\frac {\beta
\mbox{} \left( 1-2\,iv,2\,m+1 \right) }{\Gamma \left( 2\,m+1 \right) }}=\frac{1}{\Gamma \left( 2\,m+1 \right)}{ {\int_{0}^{1}\!{x}^{-2\,iv}
\mbox{} \left( 1-x \right) ^{2\,m}\,{\rm d}x}{ }}.  
\label{betaDef}
\end{equation} 

Substitute \eqref{betaDef} into \eqref{Cp4} and again reorder the integral operations to obtain

\begin{equation}
\displaystyle A(m)= {\frac {{2}^{2\,m+1}}{{\pi}^{2}\Gamma \left( 2\,m+1 \right) }\int_{0}^{\infty }\!{\frac {{t}^{2\,m+1}}{\cosh \left( t \right) }\int_{0}^{1}\! \left( 1-x \right) ^{2\,m}\int_{-\infty }^{\infty }\!{\frac {1}{\cosh \left( \pi\,v \right) } \left( {\frac {\pi}{2\,xt}} \right) ^{2\,iv}}\,{\rm d}v
\mbox{}\,{\rm d}x}\,{\rm d}t}\,.
\label{Cp4b}
\end{equation}


The innermost integration in \eqref{Cp4b} can be found from \cite[Eq. 3.981(3)]{G&R} by converting the integral into the range ($0,\infty $), writing $\left(\frac {\pi}{2xt} \right) ^{2\,iv}$ in the form of complex trigonometric functions and applying symmetry to the result, eventually giving  
\begin{equation}
\displaystyle \int_{-\infty }^{\infty }\!{\frac {1}{\cosh \left( \pi\,v \right) } \left( {\frac {\pi}{2xt}} \right) ^{2\,iv}}\,{\rm d}v
\mbox{}={\frac {4\pi\,xt}{4\,{x}^{2}{t}^{2}+{\pi}^{2}}}\,.
\label{Cint1}
\end{equation}

Thus, after rescaling the $t-$integration by $\pi$, \eqref{Cp4b} becomes

\begin{equation}
\displaystyle A \left( m \right)= \, {\frac {{2}^{2\,m+3}
\mbox{}{\pi}^{2\,m}}{\Gamma \left( 2\,m+1 \right) }\int_{0}^{\infty }\!{\frac {{t}^{2\,m+2}}{\cosh \left( \pi\,t \right) }\int_{0}^{1}\!{\frac { \left( 1-x \right) ^{2\,m}x}{4\,{t}^{2}{x}^{2}+1}}\,{\rm d}x}\,{\rm d}t}\,.
\label{Cp5}
\end{equation}

\subsubsection{Reduction of \eqref{Cp5}}

The inner integral in \eqref{Cp5} is listed in \cite[Eq. 3.254(1)]{G&R}, giving


\begin{equation}
\displaystyle \int_{0}^{1}\!{\frac { \left( 1-x \right) ^{2\,m}x}{4\,{t}^{2}{x}^{2}+1}}\,{\rm d}x={\frac {\Gamma \left( 2\,m+1 \right) {\mbox{$_3$F$_2$}(1,1,3/2;\,2+m,3/2+m;\,-4\,{t}^{2})}
\mbox{}}{\Gamma \left( 2\,m+3 \right) }}\,,
\label{Ix1}
\end{equation}
in which case, we have


\begin{equation}
\displaystyle A \left( m \right) ={\frac {{2}^{3+2\,m}
\mbox{}{\pi}^{2\,m}}{\Gamma \left( 3+2\,m \right) }\int_{0}^{\infty }\!{\frac {{t}^{2\,m+2}{\mbox{$_3$F$_2$}(1,1,3/2;\,2+m,3/2+m;\,-4\,{t}^{2})}}{\cosh \left( \pi\,t \right) }}\,{\rm d}t}\,.
\label{A(m)}
\end{equation}

This result can be verified by direct integration for at least the cases $m=1$ and $m=2$. See \cite[Section 6]{milgram2020additions}.

\subsubsection{Alternate reduction of \eqref{Cp5} and an asymptotic divergent series for Hurwitz' alternating Zeta function}


Since both integrals composing \eqref{Cp5} are convergent at both endpoints with $m>0$, rather than directly evaluating \eqref{Cp5} as in the previous section, interchange the order of integration giving
\begin{equation}
\displaystyle A \left( m \right) ={\frac {{2}^{2\,m+3}{\pi}^{2\,m}}{\Gamma \left( 2\,m+1 \right) }\int_{0}^{1}\! \left( 1-x \right) ^{2\,m}x\int_{0}^{\infty }\!{\frac {{t}^{2\,m+2}}{ \left( 4\,{t}^{2}{x}^{2}+1 \right) \cosh \left( \pi\,t \right) }}\,{\rm d}t
\mbox{}\,{\rm d}x}\,.
\label{Cp5b}
\end{equation}
From \cite[Eq. (3.10)]{milgram2020additions} with $a=1/x$, the inner integral can be evaluated to obtain:
\begin{align} \nonumber
\displaystyle \int_{0}^{\infty }\!{\frac {{t}^{2\,m+2}}{ \left( 4\,{t}^{2}+1/{x}^{2} \right) \cosh \left( \pi\,t \right) }}\,{\rm d}t&=\frac { \left( -1 \right) ^{m} }{{2}^{2\,m+3}x^{2m}} 
\left[\Gamma \left( 2\,m+2 \right) \sum _{j=0}^{2\,m}{\frac { \left( -1 \right) ^{j} \left( 2\,x \right) ^{j}\mathcal{E} \left( j,\frac{1}{2x}+\frac{1}{2} \right)}{\Gamma \left( j+2 \right) 
\mbox{}\Gamma \left( 2\,m-j+1 \right) } } 
\mbox{}\right.\\
&\left. -{\frac {1}{2\,x} \left( \psi \left(\frac{1}{4x} +\frac{3}{4}\right) -\psi \left( \frac{1}{4x}+\frac{1}{4} \right)  \right) }
\mbox{} \right]\,.
\label{Cor2a}
\end{align}

Substitute \eqref{EuIdhalf} into \eqref{Cor2a}, apply the Lemma \eqref{Lemma1}, substitute into \eqref{Cp5b} with a change of variables $x\rightarrow 1/t$ and, noting that $\mathcal{E}(2k+1)=0$, eventually obtain

\begin{equation}
\displaystyle A \left( m \right) \, = \,-\left( -1 \right) ^{m}\frac{{\pi}^{2\,m}}{2\,\Gamma \left( 2\,m+1 \right)}\,{ { \int_{1}^{\infty }\! \left( t-1 \right) ^{2\,m} \left[ \psi \left( \frac{t}{4}+\frac{3}{4} \right) -\psi \left( \frac{t}{4}+\frac{1}{4} \right) -2\,\sum _{k=0}^{m}\frac{\mathcal{E} \left( 2\,k \right)}{{t}^{2\,k+1}}  \right] 
\mbox{}\,{\rm d}t}{ }}\,.
\label{Am2}
\end{equation}

The bracketed term in \eqref{Am2} must asymptotically decrease faster than $O(t^{-2m-1})$ since the original integral was convergent at $x=0$, and nothing has happened to alter that state of affairs; thus the transformed integrand must be of $O(t^{-1-\epsilon})$ at its upper limit. From this we infer that, for $t\rightarrow \infty\,$,

\begin{equation}
\displaystyle \psi \left( \frac{t}{4}+\frac{3}{4} \right) -\psi \left( \frac{t}{4}+\frac{1}{4} \right)\sim 2\,\sum _{k=0}^{m}{\frac {\mathcal{E} \left( 2\,k \right) }{{t}^{2\,k+1}}}
\label{C1Asy}
\end{equation}
and since the left-hand side of \eqref{C1Asy} is independent of $m$, this must be true for all values of $m>0$, in which case \eqref{C1Asy} is equivalent to a divergent asymptotic series representation of the left-hand side, that is

\begin{equation}
\displaystyle \psi \left( \frac{t}{4}+\frac{3}{4} \right) -\psi \left( \frac{t}{4}+\frac{1}{4} \right)\sim 2\,\sum _{k=0}^{\infty}{\frac {\mathcal{E} \left( 2\,k \right) }{{t}^{2\,k+1}}}\,.
\label{C2Asy}
\end{equation}
This result that can be easily verified numerically by expanding the left-hand side asymptotically using either the Mathematica \cite{Math} or Maple \cite{Maple} ``Series" command, to any order desired. Note the identity (\cite[Eq. 22(5)]{Sriv&Choi})
\begin{equation}
\displaystyle \psi \left( \frac{t}{4}+\frac{3}{4} \right) -\psi \left( \frac{t}{4}+\frac{1}{4} \right)=2\eta(1,(t+1)/2)
\label{PdiId}
\end{equation}
where $\eta(1,(t+1)/2)$ is the alternating Hurwitz Zeta function (a.k.a the Lerch transcendent - $\Phi(1/2,s,a)$ - see \cite{Williams-Zhang}). Thus \eqref{C2Asy} and \eqref{PdiId} together give an infinite (asymptotic) divergent series representation of a special case of the alternating Hurwitz Zeta function to all inverse orders of $t$ - see \eqref{AsyInv}. Also see \cite[Theorem 1a]{Borwein^2Dilcher}.

\section{Derivative}\label{sec:Derivative}


The derivative  $ \frac {\partial}{\partial s} E_{s}(z) $ defines the ``Generalized Integro-Exponential function" as follows

\begin{equation}
E_{s}^{j}(z)\equiv \frac{(-1)^j}{\Gamma(j+1)}\frac {\partial^{j}}{\partial s^{j}}E_{s}(z)\,\hspace{20pt}j\ge0.
\label{EsDef}
\end{equation}

Two important properties \cite[Eqs. (2.4) and (2.8)]{Milgram:1985}  are the recursion

\begin{equation}
\displaystyle {E^j_s} \left(z \right) ={\frac {z{E^j_{s-1}} \left( z \right) -{E^{j-1}_s} \left(z \right) }{1-s}}\,,\hspace{20pt}s\neq1,j\geq0,
\label{EiGenRecur}
\end{equation}
and the limiting case
\begin{equation}
{E^{-1}_s} \left(z \right) =\exp(-z)\,.
\label{Ejm1}
\end{equation}
Thus after straightforward differentiation of \eqref{EiSum} we have
\begin{equation}
\displaystyle {\frac {\rm d}{{\rm d}s}}\eta \left( s \right) =\ln  \left( 2 \right) \eta \left( s \right) +{2}^{s-1}\sum _{k=0}^{\infty }\left({E^{1}_{s}} \left(-i\pi\, \left( k+1/2 \right)  \right) +{E^{1}_{s}} \left( i\pi\, \left( k+1/2 \right)  \right)\right) 
\label{DiffEta}
\end{equation}

which can alternatively be written as 

\begin{equation}
\displaystyle  \left( 1-{2}^{1-s} \right)\zeta^{\prime} \left( s \right) ={ {{ }\ln  \left( 2 \right)  \left(1-{2}^{2-s} \right) \zeta \left( s \right) }}
\mbox{}+{ {{2}^{s-1}\sum _{k=0}^{\infty }\left({E^{1}_{s}} \left(-i\pi\, \left( k+1/2 \right)  \right) +{E^{1}_{s}} \left( i\pi\, \left( k+1/2 \right)  \right)\right) }{}}\,.
\label{DiffZeta}
\end{equation}

Again using forward recursion, from \eqref{R2a}, \eqref{EiGenRecur} and \eqref{DiffZeta} we find the equivalent form of \eqref{R2a}, that is
\begin{align} \nonumber
\displaystyle  \left( 1-{2}^{1-s} \right) \zeta^{\prime} \left( s \right)& ={2}^{s}{\frac {is}{2\pi}\sum _{k=0}^{\infty }{\frac {{E^{1}_{s+1}} \left( i\pi\, \left( k+1/2 \right)  \right) -{E^{1}_{s+1}} \left(- i\pi\, \left( k+1/2 \right)  \right) 
\mbox{}}{k+1/2}}} \\
&+ \left( \ln  \left( 2 \right)  \left( 1-{2}^{2-s} \right) +{\frac {1-{2}^{1-s}}{s}} \right) \zeta \left( s \right)
\mbox{}-{\frac {{2}^{s-1}}{s}}\,.
\label{Rc1bD}
\end{align}
For the record, the general form corresponding to \eqref{R5} with $n\geq 0$ is

\begin{align} \nonumber
\displaystyle  \left( 1-{2}^{1-s} \right) \zeta^{\prime}& \left(s \right) ={\frac {{2}^{s-1}}{\Gamma \left( s \right) }\sum _{j=0}^{\left\lfloor n/2-1/2\right\rfloor}{\frac { \left( \psi \left( 2\,j+s \right) -\psi \left( s+n \right)  \right) 
\mbox{}\Gamma \left( 2\,j+s \right) \mathcal{E} \left( 2\,j \right) }{\Gamma \left( 2\,j+1 \right) }}}
\\ \nonumber
&+{\frac {{2}^{s-1}{{\rm e}^{i\pi n/2}} \left( -1 \right) ^{n}\Gamma \left( s+n \right) }{\Gamma \left( s \right) 
{\pi}^{n}}\sum _{k=0}^{\infty }{\frac {{E^{1}_{s+n}} \left( -i\pi\, \left( k+1/2 \right)  \right) + \left( -1 \right) ^{n}{\it E^1_{s+n}} \left(i\pi\, \left( k+1/2 \right)  \right) 
\mbox{}}{ \left( k+1/2 \right) ^{n}}}}\\&+ \left(  \left( 1-{2}^{2-s} \right) \ln  \left( 2 \right) + \left(\psi \left( s+n \right) -\psi \left( s \right)   \right)  \left( 1-{2}^{1-s} \right)  \right) \zeta \left( s \right)\,. 
\label{Rgen}
\end{align} 

\section{An Integro-series representation} \label{sec:Integro}
The functions $E_s( z)$ appearing in \eqref{ExpInt} are related to incomplete (upper) Gamma functions via 
\begin{equation}
E_s(z) = z^{-1+s}\Gamma(1-s, z)\,.
\label{EiGamma}
\end{equation}

With reference to \eqref{EiGamma}, a useful integral representation \cite[page 42]{Olver} is
\begin{equation}
\displaystyle \Gamma \left( s,\pm i\kappa \right) ={{\rm e}^{\pm is\pi/2}}\int_{\kappa}^{\infty }\!{{\rm e}^{\mp iv}}{v}^{s-1}\,{\rm d}v
\label{Olver}
\end{equation}
or, with $\kappa >0$, and a change of integration variables
\begin{equation}
\displaystyle \Gamma \left( s,\pm i\kappa \right) =\Gamma \left( s,\kappa \right) +{\kappa}^{s}\int_{\pm i}^{1}\!{{\rm e}^{-\kappa\,v}}{v}^{s-1}\,{\rm d}v\,.
\label{GamId}
\end{equation}


Keeping \eqref{EiSum} in mind, we have
\begin{align} \nonumber
\displaystyle  {{\rm e}^{is\pi/2}}\Gamma \left( s,-i\kappa \right) 
\mbox{}+{{\rm e}^{-is\pi/2}}\Gamma \left( s,i\kappa \right) &=
{{\rm e}^{is\pi/2}}{\kappa}^{s}\int_{-i}^{1}\!{{\rm e}^{-\kappa\,v}}{v}^{s-1}\,{\rm d}v+{{\rm e}^{-is\pi/2}}{\kappa}^{s}\int_{i}^{1}\!{{\rm e}^{-\kappa\,v}}{v}^{s-1}\,{\rm d}v \\&+ 2\,\cos \left(s\pi/2 \right) \Gamma \left( s,\kappa \right) \,,
\label{Eq2}
\end{align}
which, with \eqref{EiGamma}, gives
\begin{align} 
\displaystyle {\it E_s} \left(i \kappa  \right) +{\it E_s} \left(-i \kappa  \right) =i{{\rm e}^{-is\pi/2}}\,\int_{-i}^{1}\!{{\rm e}^{-\kappa v}}{v}^{-s}\,{\rm d}v
\mbox{}-i{{\rm e}^{is\pi/2}}\,\int_{i}^{1}\!{{\rm e}^{-\kappa v}}{v}^{-s}\,{\rm d}v+2 \sin \left( s\pi/2 \right){\it E_s} \left(\kappa  \right) 
\label{Eq4}
\end{align}
after setting $s\rightarrow 1-s$. Further, the identification $\kappa=\pi(k+1/2)$, interchange of integration and summation (reversing the previous interchange), together with the simple identity
\begin{equation}
\displaystyle \sum _{k=0}^{\infty }{{\rm e}^{-\pi\, \left( k+1/2 \right) v}}=\frac{1}{2\, \sinh \left( \pi\,v/2 \right)  } 
\label{SimpleSum}
\end{equation}
gives
\begin{align} \nonumber
\displaystyle \sum _{k=0}^{\infty }&{E_s} \left(i\pi\, \left( k+1/2 \right)  \right) +{ E_s} \left( -i\pi\, \left( k+1/2 \right)  \right) 
\mbox{}=\frac{i{{\rm e}^{-is\pi/2}}}{2}\int_{-i}^{1}\!{\frac {{v}^{-s}}{\sinh \left( \pi\,v /2\right) }}\,{\rm d}v
\mbox{}\\
&-\frac{i{{\rm e}^{is\pi/2}}}{2}\int_{i}^{1}\!{\frac {{v}^{-s}}{\sinh \left(\pi\,v/2 \right) }}\,{\rm d}v+2\,\sin \left( s\pi/2 \right) \sum _{k=0}^{\infty }{E_{s}} \left(\pi\, \left( k+1/2 \right)  \right)\,. 
\label{BasicEq}
\end{align}
Thus, from \eqref{EiSum} and \eqref{BasicEq} we have

\begin{align} \nonumber
\displaystyle  \frac {{2}^{1-s}-1}
{{2}^{s-1}} \zeta \left( s \right)& =  \frac{i{{\rm e}^{-is\pi/2}}}{2}\int_{-i}^{1}\,{\frac {{v}^{-s}}{\sinh \left( \pi\,v /2\right) }}\,{\rm d}v
+\frac{i{{\rm e}^{is\pi/2}}}{2}\int_{1}^{i}\,{\frac {{v}^{-s}}{\,\sinh \left(\pi\,v/2 \right) }}\,{\rm d}v
\\&+2\,\sin \left( s\pi/2 \right) \sum _{k=0}^{\infty }{E_{s}} \left(\pi\, \left( k+1/2 \right)  \right)\,.
\label{ZetaSum}   
\end{align}
Notice that the structure of $\zeta(s)$ is fundamentally different for $s=2n$ or $s=2n+1$ because of the coefficient $\sin(\pi s/2)$ in \eqref{ZetaSum}. At this point, there are a number of interesting possible choices for the integration paths in \eqref{ZetaSum}, each of which leads to different identities. See Figure \ref{fig:Figure1}. \newline

\begin{figure}[h] 
\centering
\includegraphics[width=.6\textwidth]{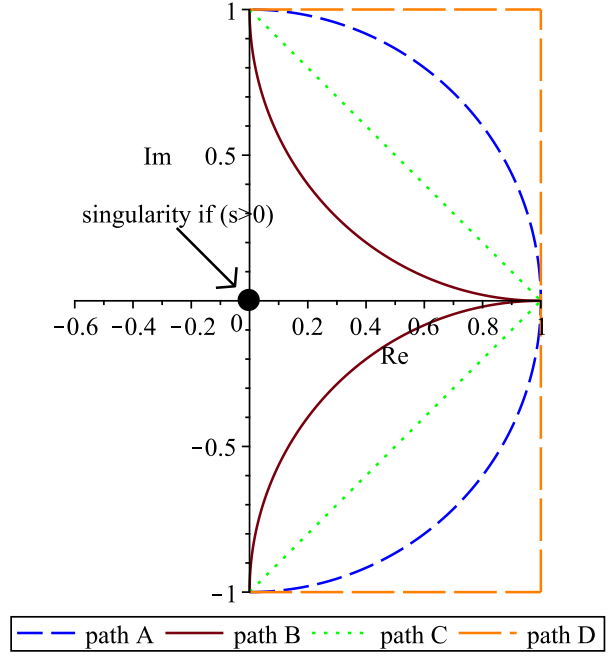}
\caption{Four possible paths of integration from $-i$ to $i$. }
\label{fig:Figure1}
\end{figure}

\subsection{Integration path A: $v={\rm e}^{i\theta}$} \label{sec:circle}

A simple choice suggests a change of variables $v={\rm e}^{i\theta}$ with $-\pi/2< \theta<\pi/2$, equivalent to choosing a circular path in the right-hand complex $v$-plane centred at the origin passing through the limiting points of the integrals in \eqref{BasicEq} (see Figure \ref{fig:Figure1}). With this choice, and following simplification (specifically, reverse one integral (i.e. $\theta\rightarrow \pi/2-\theta$) and utilize symmetry under $\theta\rightarrow -\theta$), \eqref{EiSum} becomes

\begin{align} \nonumber
\displaystyle {\frac {\zeta \left( s \right)  \left( 1-{2}^{1-s} \right) }{{2}^{s}}}&=2^{-s}\eta(s)=-\sin \left( s\pi/2 \right) \sum _{k=0}^{\infty }{\it E_s} \left(\pi \left( k+1/2 \right)  \right)\\
&+\frac{1}{2}\,\int_{-\pi/2}^{\pi/2}\!{\frac {{{\rm e}^{-(\pi/2)\,\sin \left( \theta \right) }}\sin \left( (\pi/2)\,\cos \left( \theta \right) +\theta\, \left( s-1 \right)  \right) 
\mbox{}}{\cosh \left( \pi\,\sin \left( \theta \right)  \right) -\cos \left( \pi\,\cos \left( \theta \right)  \right) }}\,{\rm d}\theta \,.
\label{Sa2b}
\end{align}

Also, see \cite[Ex. 3.2, page 43]{Olver}. From Appendix A, we thus have another representation for the lower incomplete function $\zeta^{-}(s)$ (see \eqref{ZetaMinus}).

\begin{equation}
\displaystyle \zeta _{{-}} \left( s \right)={2}^{s-1}\int_{-\pi/2}^{\pi/2}\!{\frac { {{\rm e}^{-(\pi/2)\,\sin \left( \theta \right) }}\sin \left( (\pi/2)\,\cos \left( \theta \right) +\theta\, \left( s-1 \right)  \right)
\mbox{}}{\cosh \left( \pi\,\sin \left( \theta \right)  \right) -\cos \left( \pi\,\cos \left( \theta \right)  \right) }}\,{\rm d}\theta
\label{ZmEq1}
\end{equation}
valid for all $s$.\newline

{\bf Corollary.}\newline In \eqref{Sa2b}, replace $s\rightarrow s+1$ and concurrently replace $s\rightarrow 1-s$, add the two results and simplify to eventually obtain
\begin{align} 
\displaystyle \int_{-\pi/2}^{\pi/2}\!{\frac {{{\rm e}^{-(\pi/2)\,\sin \left( \theta \right) }}\sin \left( (\pi/2)\cos \left( \theta \right)  \right) \cos \left( s\,\theta \right) }{\cosh \left( \pi\,\sin \left( \theta \right)  \right) 
-\cos \left( \pi\,\cos \left( \theta \right)  \right) }}{\rm d}\theta&={2}^{-s-1} \left( 1-{2}^{-s} \right) {\zeta^{-}} \left( 1+s \right)+{2}^{s-1} \left( 1-{2}^{s} \right) { \zeta^{-}} \left( 1-s \right) \\
&= 2^{-s-1}\eta^{-}(1+s)+2^{s-1}\eta^{-}(1-s)\,.
\label{Sa2f}
\end{align}

Similarly, by subtracting, we identify

\begin{align}
\displaystyle \int_{-\pi/2}^{\pi/2}\!{\frac {{{\rm e}^{-(\pi/2)\sin \left( \theta \right) }}\cos \left( (\pi/2)\cos \left( \theta \right)  \right) 
\mbox{}\sin \left( s\,\theta \right) }{\cosh \left( \pi\,\sin \left( \theta \right)  \right) -\cos \left( \pi\,\cos \left( \theta \right)  \right) }}{\rm d}\theta&={2}^{-s-1} \left( 1-{2}^{-s} \right) {\zeta^{-}} \left( 1+s \right)-{2}^{s-1} \left( 1-{2}^{s} \right) { \zeta^{-}} \left( 1-s \right)\\
&={2}^{-s-1}{\it \eta^{-}} \left( 1+s \right) -{2}^{s-1}{\it \eta^{-}} \left( 1-s \right) 
\label{Sa2h}
\end{align}

{\bf Remark:} The following subsections list some special cases that may be worthy of special attention. It is obvious that variations of all of the following can be obtained by simple transformations of integration variables (e.g. $\theta\rightarrow \pi/2-\theta$ etc.)
\subsubsection{Special case: s=1}

In the case that $s=1$, take the limit of the left-hand side of \eqref{Sa2b} and eventually obtain
\begin{equation}
\displaystyle \int_{-\pi/2}^{\pi/2}\!{\frac {\sin \left( (\pi/2)\,\cos \left( \theta \right)  \right) {{\rm e}^{-(\pi/2)\,\sin \left( \theta \right) }}}{\cosh \left( \pi\sin \left( \theta \right)  \right) 
\mbox{}-\cos \left( \pi\,\cos \left( \theta \right)  \right) }}\,{\rm d}\theta=\ln  \left( 2 \right) +2\,\sum _{k=0}^{\infty } E_{1} \left(\pi \left( k+1/2 \right)  \right) 
\label{Scase1a}
\end{equation}
which can also be written
\begin{align}
\displaystyle \int_{0}^{\pi/2}\!{\frac {\sin \left( (\pi/2)\,\cos \left( \theta \right)  \right) \cosh \left( (\pi/2)\,\sin \left( \theta \right)  \right) }{\cosh \left( \pi\,\sin \left( \theta \right)  \right) 
\mbox{}-\cos \left( \pi\,\cos \left( \theta \right)  \right) }}\,{\rm d}\theta=\frac{1}{2}\ln  \left( 2 \right) +\sum _{k=0}^{\infty }{ E_{1}} \left(\pi\, \left( k+1/2 \right)  \right). 
\label{Scase1}
\end{align}
The sum in both of these can be identified from \eqref{Zplus} with $s=1$, noting that $\eta(1)$ does not share the singularity of $\zeta(1)$.

\subsubsection{Special cases $s=\pm 2n$ and $s=1\pm 2n$}


The case $s=0,\pm 1,\pm 2\dots$ leads to some interesting identifications from the well-known identification (\cite[Eq. 25.6.2]{NIST}) of $\zeta(2n)$ in terms of Bernoulli numbers $B_{2n}$ .\newline

If $0\leq s=2n$ then
\begin{equation}
\displaystyle \int_{-\pi/2}^{\pi/2}\!{\frac { {{\rm e}^{-(\pi/2)\sin \left( \theta \right) }}\sin \left( (\pi/2)\cos \left( \theta \right) 
+(2\,n-1)\,\theta \right)
\mbox{}}{\cosh \left( \pi\,\sin \left( \theta \right)  \right) -\cos \left( \pi\,\cos \left( \theta \right)  \right) }}\,{\rm d}\theta=\left( -1 \right) ^{n}{B_{2n\, }{\pi}^{2\,n}\frac {  
\mbox{} \left( {2}^{1-2\,n}-1 \right) }{\Gamma \left( 2\,n+1 \right) }}\,.
\label{Sa2cEp}
\end{equation}
In the above, if $n=0$, we have
\begin{equation}
\displaystyle \int_{-\pi/2}^{\pi/2}\!{\frac {{{\rm e}^{-(\pi/2)\sin \left( \theta \right) }}\sin \left( (\pi/2)\cos \left( \theta \right) -\theta \right) }{\cosh \left( \pi\,\sin \left( \theta \right)  \right) 
\mbox{}-\cos \left( \pi\,\cos \left( \theta \right)  \right) }}\,{\rm d}\theta=1\,.
\label{Casen0}
\end{equation}

If $s=-2n<0$, then 
\begin{equation}
\displaystyle \int_{-\pi/2}^{\pi/2}\!{\frac {{{\rm e}^{-(\pi/2)\sin \left( \theta \right) }}\,\sin \left( (\pi/2)\cos \left( \theta \right) 
-(2\,n+1)\,\theta \right) 
\mbox{}}{\cosh \left( \pi\,\sin \left( \theta \right)  \right) -\cos \left( \pi\,\cos \left( \theta \right)  \right) }}\,{\rm d}\theta=0\,.
\label{Sa2cEm}
\end{equation}
If $1\leq s=2n+1$ then
\begin{align} \nonumber
\displaystyle \int_{-\pi/2}^{\pi/2}\!{\frac {{{\rm e}^{-(\pi/2)\sin \left( \theta \right) }}\,\sin \left( (\pi/2)\cos \left( \theta \right) +2\,n\,\theta \right) }{\cosh \left( \pi\,\sin \left( \theta \right)  \right) 
\mbox{}-\cos \left( \pi\,\cos \left( \theta \right)  \right) }}\,{\rm d}\theta&={4}^{-n} \left( 1-{4}^{-n} \right) \zeta \left( 2\,n+1 \right) \\
&+2\, \left( -1 \right) ^{n}\sum _{k=0}^{\infty }{ E_{ 2n+1}} \left(\pi\, \left( k+1/2 \right)  \right) \,,
\label{Sa2cOp}
\end{align} 
and, in the above, if $n=0$, by taking a limit, we obtain
\begin{equation}
\displaystyle \int_{-\pi/2}^{\pi/2}\!{\frac {{{\rm e}^{-(\pi/2)\sin \left( \theta \right) }}\sin \left( (\pi/2)\cos \left( \theta \right)  \right) }{\cosh \left( \pi\,\sin \left( \theta \right)  \right) 
\mbox{}-\cos \left( \pi\,\cos \left( \theta \right)  \right) }}\,{\rm d}\theta=\ln  \left( 2 \right) +2\,\sum _{k=0}^{\infty }{\it E_{1}} \left(\pi\, \left( k+1/2 \right)  \right) \,.
\label{Sa2cOpA}
\end{equation}
Finally, if $s=-2n+1<0$ then
\begin{align} \nonumber
\displaystyle \int_{-\pi/2}^{\pi/2}\!{\frac {{{\rm e}^{-(\pi/2)\sin \left( \theta \right) }}\,\sin \left( (\pi/2)\cos \left( \theta \right) -2\,\theta\,n \right) 
\mbox{}}{\cosh \left( \pi\,\sin \left( \theta \right)  \right) -\cos \left( \pi\,\cos \left( \theta \right)  \right) }}\,{\rm d}\theta&={2}^{2\,n-1}B_{2n\,}{\frac {  \left({2}^{2\,n}  -1\right) }{n}}
\\
&+2\, \left( -1 \right) ^{n}\sum _{k=0}^{\infty }{ E_{1-2\,n}} \left( \pi\, \left( k+1/2 \right)  \right) \,.
\label{Sa2cOm}
\end{align}
In \eqref{Sa2cOm}, recall that for $m=1,2\dots$
\begin{equation}
\displaystyle {\it E_{-m}} \left(\kappa \right) ={\frac {\Gamma \left( m+1 \right) {{\rm e}^{-\kappa}}}{{\kappa}^{m+1}}\sum _{j=0}^{m}{\frac {{\kappa}^{j}}{j!}}}=\Gamma(1+m,\kappa)/\kappa^{m+1}\,,
\label{EiNeg_m}
\end{equation}
which, together with \eqref{EiPoly} and \eqref{LiLn} yields an alternate form for \eqref{Sa2cOm}

\begin{align} \nonumber
&\displaystyle \int_{-\pi/2}^{\pi/2}\!{\frac {{{\rm e}^{-(\pi/2)\sin \left( \theta \right) }}\sin \left( (\pi/2)\cos \left( \theta \right) -2n\,\theta\right) }{\cosh \left( \pi\,\sin \left( \theta \right)  \right) 
\mbox{}-\cos \left( \pi\,\cos \left( \theta \right)  \right) }}\,{\rm d}\theta= {2}^{2\,n} \left( 1-{2}^{2\,n} \right)\,\zeta \left( 1-2\,n \right) \\
&+{\frac { 2\,\left( -1 \right) ^{n}\Gamma \left( 2\,n \right) }{\pi}\sum _{j=1}^{2\,n-1}{\frac {{2}^{j+1}\,Li_{j+1} \left({{\rm e}^{-\pi/2}} \right) -Li_{j+1} \left({{\rm e}^{-\pi}} \right) }{\Gamma \left( 2\,n-j \right) 
\mbox{}{\pi}^{j}}}}+{\frac {2\, \left( -1 \right) ^{n} }{\pi}}\ln  \left( \displaystyle {\frac {\sinh \left( \pi/2 \right) }{\cosh \left( \pi/2 \right) -1}} \right)\,.
\label{Eq4A}
\end{align}

Adding and subtracting then yields the following identities
\begin{align} \nonumber
\displaystyle  &\int_{-\pi/2}^{\pi/2}\!{\frac {{{\rm e}^{-(\pi/2)\sin \left( \theta \right) }}\cos \left( 2n\,\theta \right) \sin \left( (\pi/2)\cos \left( \theta \right) -\theta \right) 
}{\cosh \left( \pi\,\sin \left( \theta \right)  \right) -\cos \left( \pi\,\cos \left( \theta \right)  \right) }}\,{\rm d}\theta\\
&=\int_{-\pi/2}^{\pi/2}\!{\frac {{{\rm e}^{-(\pi/2)\sin \left( \theta \right) }}\sin \left( 2n\,\theta \right) \cos \left( (\pi/2)\cos \left( \theta \right) -\theta \right) }{\cosh \left( \pi\,\sin \left( \theta \right)  \right) 
\mbox{}-\cos \left( \pi\,\cos \left( \theta \right)  \right) }}\,{\rm d}\theta  = \left( -1 \right) ^{n} B_{  2n }{\pi}^{2\,n}  {\frac {  \left( {2}^{1-2\,n}-1 \right) 
\mbox{}}{2\Gamma \left( 2\,n+1 \right) }} \,, 
\label{Sa2Epm}
\end{align}
\begin{align} \nonumber
\displaystyle & \int_{-\pi/2}^{\pi/2}\!{\frac {{{\rm e}^{-(\pi/2)\sin \left( \theta \right) }}\cos \left( 2n\,\theta \right) \sin \left( (\pi/2)\cos \left( \theta \right)  \right) 
\mbox{}}{\cosh \left( \pi\,\sin \left( \theta \right)  \right) -\cos \left( \pi\,\cos \left( \theta \right)  \right) }}\,{\rm d}\theta
=
\zeta \left( 2\,n+1 \right) {2}^{-2\,n-1} \left( 1-{4}^{-n} \right)\\
&+ \left( -1 \right) ^{n}\sum _{k=0}^{\infty }{E_{ 2n+1}} \left(\pi \left( k+1/2 \right)  \right) +
{ E_{1-2\,n}} \left( \pi\, \left( k+1/2 \right)  \right)+{{B_{2n\,}}\frac {{2}^{2\,n-2} \left({2}^{2\,n}-1 \right)   }{n}} \,,
\label{Sa2Op}
\end{align}
and
\begin{align} \nonumber
\displaystyle &\int_{-\pi/2}^{\pi/2}\!{\frac {{{\rm e}^{-(\pi/2)\sin \left( \theta \right) }}\sin \left( 2n\,\theta \right) \cos \left( (\pi/2)\cos \left( \theta \right)  \right) 
\mbox{}}{\cosh \left( \pi\,\sin \left( \theta \right)  \right) -\cos \left( \pi\,\cos \left( \theta \right)  \right) }}\,{\rm d}\theta=\zeta \left( 2\,n+1 \right) {2}^{-2\,n-1} \left( 1-{4}^{-n} \right) \\
&+\left( -1 \right) ^{n}\sum _{k=0}^{\infty } { E_{ 2n+1}} \left(\pi \left( k+1/2 \right)  \right) -E_{1-2n} \left(\pi \left( k+1/2 \right)  \right) 
\mbox{}-{B_{2n}\frac {{2}^{2\,n-2} \left( {2}^{2\,n}-1 \right)   }{n}}\,.
\label{Sa2Om}
\end{align}

\subsubsection{Special cases s=0, s=1 and s=2 mixed}

By adding and subtracting cases corresponding to $s=0$ and $s=1$, after some simplification we obtain

\begin{equation}
\displaystyle \int_{-\pi/2}^{\pi/2}\!{\frac {{{\rm e}^{-(\pi/2)\sin \left( \theta \right) }}\sin \left( \theta/2 \right) \cos \left( (\pi/2)\cos \left( \theta \right) -\theta/2 \right) 
\mbox{}}{\cosh \left( \pi\,\sin \left( \theta \right)  \right) -\cos \left( \pi\,\cos \left( \theta \right)  \right) }}\,{\rm d}\theta=\ln  \left( 2 \right)/2 -1/2+\sum _{k=0}^{\infty }{ E_{1}} \left(\pi\, \left( k+1/2 \right)  \right)\,.
\label{Mixeds01}
\end{equation}

Similarly, adding and subtracting the cases s=0 and s=2 along with some elementary trigonometric identities, respectively yields

\begin{equation}
\displaystyle \int_{-\pi/2}^{\pi/2}\!{\frac {\cos \left( \theta \right) \sin \left( (\pi/2)\,\cos \left( \theta \right)  \right) {{\rm e}^{-(\pi/2)\,\sin \left( \theta \right) }}}{\cosh \left( \pi\,\sin \left( \theta \right)  \right) 
\mbox{}-\cos \left( \pi\,\cos \left( \theta \right)  \right) }}\,{\rm d}\theta=1/2+{\pi}^{2}/48\,,
\label{V02p}
\end{equation}
and
\begin{equation}
\displaystyle \int_{-\pi/2}^{\pi/2}\!{\frac {\sin \left( \theta \right) \cos \left( (\pi/2)\,\cos \left( \theta \right)  \right) {{\rm e}^{-(\pi/2)\sin \left( \theta \right) }}}{\cosh \left( \pi\,\sin \left( \theta \right)  \right) 
\mbox{}-\cos \left( \pi\,\cos \left( \theta \right)  \right) }}\,{\rm d}\theta=-1/2+{\pi}^{2}/48\,.
\label{V02m}
\end{equation}

\subsubsection{Special case s=1/2}


The case $s=1/2$ in \eqref{Sa2b} immediately yields the identity 
\begin{align} \nonumber
\displaystyle \int_{-\pi/2}^{\pi/2}\!{\frac {{{\rm e}^{-(\pi/2)\sin \left( \theta \right) }}\sin \left( (\pi/2)\cos \left( \theta \right) -\theta/2 \right) }{\cosh \left( \pi\,\sin \left( \theta \right)  \right) 
\mbox{}-\cos \left( \pi\,\cos \left( \theta \right)  \right) }}\,{\rm d}\theta&= \sqrt{2} \left(  \left( 1- \sqrt{2} \right) \zeta \left( 1/2 \right) +\sum _{k=0}^{\infty }{E_{1/2}} \left( \pi\,(k+1/2) \right)  \right) \\
&=\sqrt{2}   \left( 1- \sqrt{2} \right) \zeta^{-} \left( 1/2 \right)
\label{sHalf0}
\end{align}
the second equality arising by virtue of \eqref{ZetaMinus}. An alternate representation for $\zeta(1/2)$ can be found by employing the identity (\cite[Eq. 7.11.3]{NIST})


\begin{equation}
\displaystyle { E_{1/2}} \left( z \right) = \sqrt{{\frac {\pi}{z}}} \left( 1-{\rm erf} \left( \sqrt{z}\right) \right) \,.
\label{erfc}
\end{equation}
From the integral representation of the error function (\cite[Eq. 7.2.1]{NIST}) applied  to \eqref{EiSum} after some simplification, we find

\begin{equation}
\displaystyle \zeta \left( 1/2 \right) ={\frac {1}{1- \sqrt{2}}\sum _{k=0}^{\infty } \left( 2\, \sqrt{2}\int_{0}^{1}\!\cos \left( \pi\, \left( k+1/2 \right) {t}^{2} \right) \,{\rm d}t
\mbox{}- 1/  \sqrt{k+1/2} \right)   }
\label{CosSqSum}
\end{equation}
leading to a series representation 
\begin{equation}
\displaystyle \zeta \left( 1/2 \right) ={\frac {1}{1- \sqrt{2}}\,\sum _{k=0}^{\infty }{\frac {2\,{\mathcal{F}_{c}} \left(  \sqrt{2\,k+1}\, \right) -1}{\sqrt {(k + 1/2)}}}}\,.
\label{FresnelC}
\end{equation}

where $\mathcal{F}_{c}$ is the Fresnel $C$ function (\cite[Eq. 7.2.7]{NIST}). From its asymptotic series (\cite{Maple}), it can be shown that the terms of the alternating series \eqref{FresnelC} converge as $O(k^{-1})$ and the series is therefore convergent providing that all its terms are grouped as shown.\newline

Finally, differentiating \eqref{Sa2b} with respect to $s$, setting $s=1/2$ and utilizing \eqref{EsDef}, gives

\begin{align} \nonumber
\displaystyle & \int_{-\pi/2}^{\pi/2}\!{\frac {{{\rm e}^{-(\pi/2\,)\sin \left( \theta \right) }}\cos \left( (\pi/2)\cos \left( \theta \right) -\theta/2 \right) 
\mbox{}}{\cosh \left( \pi\,\sin \left( \theta \right)  \right) -\cos \left( \pi\,\cos \left( \theta \right)  \right) }}\,\theta\,{\rm d}\theta=- \sqrt{2}\sum _{k=0}^{\infty }{E_{1/2}^{1}} \left(\pi\, \left( k+1/2 \right)  \right)\\
&+ \left(  \left( \sqrt{2}/2\, +1 \right) \ln  \left( 2 \right) + \left( \gamma+\pi/2+\ln  \left( \pi \right)  \right)  \left( \sqrt{2}/2\, -1 \right) 
\mbox{} \right) \zeta \left( 1/2 \right) -\pi\,\, \left( \sqrt{2}/2\, -1 \right){ \zeta^{+}(1/2)}
\label{ShalfDDs}
\end{align}

\subsubsection{On the critical line: s=1/2+it}

An interesting set of results can be obtained by setting $s=1/2+it$ in \eqref{Sa2b} and identifying $E_{1/2+it}(i\pi\kappa)$ by its integral representation \eqref{ExpInt}. Define

\begin{align}
&\displaystyle { J_1(t)}\, \equiv \,\int_{-\pi/2}^{\pi/2}\!{\frac {{{\rm e}^{-(\pi/2)\,\sin \left( \theta \right) }}\cos \left( (\pi/2)\,\cos \left( \theta \right) -\theta/2 \right) 
\mbox{}\sinh \left( \theta\,t \right) }{\cosh \left( \pi\,\sin \left( \theta \right)  \right) -\cos \left( \pi\,\cos \left( \theta \right)  \right) }}\,{\rm d}\theta  \label{J1}  \\ 
&\displaystyle { J_2(t)}\,\equiv \,\int_{-\pi/2}^{\pi/2}\!{\frac {{{\rm e}^{-(\pi/2)\,\sin \left( \theta \right) }}\sin \left( (\pi/2)\,\cos \left( \theta \right) -\theta/2 \right)
\mbox{}\cosh \left( \theta\,t \right) }{\cosh \left( \pi\,\sin \left( \theta \right)  \right) -\cos \left( \pi\,\cos \left( \theta \right)  \right) }}\,{\rm d}\theta \\
&\displaystyle { J_3(t)}\,\equiv\,\int_{1}^{\infty }\!{\frac {\sin \left( t\ln  \left( v \right)  \right) }{ \sqrt{v}\sinh \left( \pi\,v/2 \right) }}\,{\rm d}v\\
&\displaystyle { J_4(t)}\, \equiv\,\int_{1}^{\infty }\!{\frac {\cos \left( t\ln  \left( v \right)  \right) }{ \sqrt{v}\sinh \left( \pi\,v/2 \right) }}\,{\rm d}v\,.
\label{J4}
\end{align}

In \eqref{Sa2b}, split the resulting identities into real and imaginary components, and solve for $J_1(t)$ and $J_2(t)$ to yield 

\begin{align} \nonumber
&\displaystyle {\it J_1(t)}=H(t)\,\sin \left( t\ln  \left( 2 \right)  \right) \zeta  _{{R}}\left( 1/2+it \right)+ \left(2  -H(t)\,\cos \left( t\ln  \left( 2 \right)  \right) \right) \zeta_{{I}} \left( 1/2+it \right) 
\mbox{}\\
&\hspace{25pt}-{{\it J_3(t)}\,\frac {\cosh \left( \pi\,t/2 \right) }{ \sqrt{2}}}+{{\it J_4(t)}\,\frac {\sinh \left( \pi\,t/2 \right) }{ \sqrt{2}}}
\label{J1Soln} \\ \nonumber
&\displaystyle {\it J_2(t)}= \left(2+ H(t)\cos \left( t\ln  \left( 2 \right)  \right)  \right) \zeta_{{R}} \left( 1/2+it \right) +H(t)\sin \left( t\ln  \left( 2 \right)  \right) \zeta _{{I}} \left( 1/2+it \right)
\mbox{}\\ 
&\hspace{25pt}+{{\it J_3(t)}\,\frac {\sinh \left( \pi\,t/2 \right) }{ \sqrt{2}}}+{\it J_4(t)}\,{\frac {\cosh \left( \pi\,t/2 \right) }{ \sqrt{2}}}
\label{J2Soln}
\end{align}
where
\begin{equation}
\displaystyle H(t)\equiv\,\,4\,\cos \left( t\ln  \left( 2 \right)  \right)- \sqrt{2} \,.
\label{HConst}
\end{equation}

From the definition \eqref{ZetaPlus}, each of the integrals $J_{3}(t)$ and $J_{4}(t)$ can be connected with the real and imaginary parts of the upper incomplete function $\zeta^{+}(s)$ as follows:
\begin{align}
&\displaystyle \zeta _{{R}}^{+} \left(1/2+it \right)=\hspace{10pt}J_{{3}} \left( t \right) Q _{{1}} \left( t \right)+J_{{4}} \left( t \right) Q_{{2}} \left( t \right) 
\label{Eq1a}\\
&\displaystyle \zeta_{{I}}^{+} \left(1/2+it \right) =\,\,-J _{{3}} \left( t \right)Q_{{2}} \left( t \right) +J_{{4}} \left( t \right) Q_{{1}} \left( t \right) \,,
\label{Eq2a}
\end{align}
where
\begin{align}
&\displaystyle Q _{{1}} \left( t \right)=\,\,-P _{{1}} \left( t \right)\cosh \left( \pi\,t/2 \right)/2 -P_{{2}} \left( t \right) \sinh \left( \pi\,t/2 \right)/2 
\label{Q1Def}\\
&\displaystyle Q_{{2}} \left( t \right) =\,\,-P_{{2}} \left( t \right) \cosh \left( \pi\,t/2 \right)/2 +P_{{1}} \left( t \right) \sinh \left( \pi\,t/2 \right)/2\,, 
\label{Q2Def}
\end{align}
and
\begin{align}
&\displaystyle P_{{1}} \left( t \right) ={\frac { \left( 2\, \sqrt{2}\cos \left( t\ln  \left( 2 \right)  \right) -1 \right) \sin \left( t\ln  \left( 2 \right)  \right) 
\mbox{}}{2\, \sqrt{2}\cos \left( t\ln  \left( 2 \right)  \right) -3}}
\label{P1Def}\\
&\displaystyle P_{{2}} \left( t \right) ={\frac {2\, \sqrt{2} \left( \cos \left( t\ln  \left( 2 \right)  \right)  \right) ^{2}- \sqrt{2}-\cos \left( t\ln  \left( 2 \right)  \right) }{2\, \sqrt{2}\cos \left( t\ln  \left( 2 \right)  \right) 
\mbox{}-3}}\,.
\label{P2Def}
\end{align}
Inverting \eqref{Eq1a} and \eqref{Eq2a} identifies the integrals $J_{{3}} \left( t \right)$ and $J_{{4}} \left( t \right)$ in terms of the real and imaginary components of the upper incomplete function $\zeta^{+}(1/2+it)$ as follows:
\begin{align}
&\displaystyle J_{{3}} \left( t \right) ={\frac { 2\,\sqrt{2} \left( -4\,\cos \left( t\ln  \left( 2 \right)  \right) +3\, \sqrt{2} \right) 
\mbox{} \left( Q _{{1}} \left( t \right){\zeta_{{R}}^{+}} \left( 1/2+it \right) -Q_{{2}}\left( t \right) {\zeta_{{I}}^{+}} \left( 1/2+it \right)  \right) }{\cosh\mbox{} \left( \pi\,t \right) }}
\label{J3Soln} \\
&\displaystyle J_{{4}} \left( t \right) {\frac {=2\, \sqrt{2} \left( -4\,\cos \left( t\ln  \left( 2 \right)  \right) +3\, \sqrt{2} \right) 
\mbox{} \left( Q_{{1}} \left( t \right) {\zeta_{{I}}^{+}} \left( 1/2+it \right)+Q_{{2}} \left( t \right){ \zeta_{{R}}^{+}} \left( 1/2+it \right)  \right) }{\cosh \mbox{} \left( \pi\,t \right) }}\,.
\label{J4Soln}
\end{align}

Further simplification is possible by utilizing  \eqref{J3Soln}, \eqref{J4Soln} and \eqref{zAdd} in \eqref{J1Soln} and \eqref{J2Soln} to eliminate $J_{{3}} \left( t \right)$ and $J_{{4}} \left( t \right)$. After considerable (Maple) simplification we arrive at the identities 
\begin{align}
&\displaystyle J_{{1}} \left( t \right) =H(t)\sin \left( t\ln  \left( 2 \right)  \right) {\zeta_{{R}}^{-}}(1/2+it)+ \left( 2-H(t)\cos \left( t\ln  \left( 2 \right)  \right)  \right) {\zeta_{{I}}^{-}(1/2+it)}\,
\label{J1Solna} \\
&\displaystyle J_{{2}} \left( t \right) = \left( 2-H(t)\cos \left( t\ln  \left( 2 \right)  \right)  \right) {\zeta_{{R}}^{-}(1/2+it)}-H(t)\sin \left( t\ln  \left( 2 \right)  \right) {\zeta_{{I}}^{-}(1/2+it)}\,.
\label{J2Solna}
\end{align}
Alternatively, inversion of \eqref{J1Solna} and \eqref{J2Solna} produces an integral representation for each component of the lower incomplete function $\zeta^{-}(1/2+it)$. Specifically

\begin{align}
\displaystyle { \zeta_{{R}}^{-}}(1/2+it)=-{\frac { \left( H(t)\cos \left( t\ln  \left( 2 \right)  \right) -2 \right) J_{{2}} \left( t \right) }{6-4\, \sqrt{2}\cos \left( t\ln  \left( 2 \right)  \right) 
\mbox{}}}+{\frac {H(t)J_{{1}} \left( t \right)\sin \left( t\ln  \left( 2 \right)  \right) }{6-4\, \sqrt{2}\cos \left( t\ln  \left( 2 \right)  \right) }}
\label{ZmRSoln} \\
\displaystyle { \zeta_{{I}}^{-}}(1/2+it)=-{\frac { \left( H(t)\cos \left( t\ln  \left( 2 \right)  \right) -2 \right) J_{{1}} \left( t \right) }{6-4\, \sqrt{2}\cos \left( t\ln  \left( 2 \right)  \right) 
\mbox{}}}-{\frac {H(t)J_{{2}} \left( t \right) \sin \left( t\ln  \left( 2 \right)  \right) }{6-4\, \sqrt{2}\cos \left( t\ln  \left( 2 \right)  \right) }}\,.
\label{ZmISoln}
\end{align} 

{\bf Remark:} As opposed to the representation \eqref{ZetaMinus} that does not converge when $s=1/2$ and thus required a continuation to the representation \eqref{redefine}, setting $t=0$ reduces \eqref{ZmRSoln} to the convergent representation \eqref{sHalf0}. In other words, \eqref{ZmRSoln} and \eqref{ZmISoln} are convergent representations of  the incomplete function $\zeta^{-}(s)$ on the critical line $s=1/2+it$.
 
\subsection{Integration path B: $v=1\pm i - {\rm e}^{i\theta}$} \label{sec:InvCircle}

A second choice for the integration path in \eqref{ZetaSum} corresponds to integrating over circular paths centred at $(i,\pm i)$ for each of the two integrals in \eqref{BasicEq} respectively (see Figure \ref{fig:Figure1}). After some simplification, this yields the identity

\begin{align} \nonumber
\displaystyle &{\frac {\zeta \left( s \right)  \left( {2}^{1-s}-1 \right) }{{2}^{s-1}}}=
2\,\sin \left(s\pi/2 \right) \sum _{k=0}^{\infty }{E_{s}} \left( \pi\, \left( k+1/2 \right)  \right)\\
&+\int_{0}^{\pi/2}\!{\frac { (3-2\cos \left( \theta \right) -2\sin \left( \theta \right) )^{-s/2} \left( \cos ( \pi B/2-sA-\theta )\, {{\rm e}^{\pi C/2}}-\cos ( \pi B/2+sA+\theta )\, {{\rm e}^{-\pi C/2}}
 \right) }{\cosh \left( \pi\,C \right) -\cos \left( \pi B \right) }}{\rm d}\theta,
\label{EqAB1a}
\end{align}
where
\begin{align*}
&C=\cos(\theta)-1\\
&B=\sin(\theta)-1\\
&A=\arctan(C/B)\,.
\end{align*}
To simplify special cases, it is worth noting the following identities
\begin{align}
\displaystyle \sin \left( A \right) ={\frac {1-\cos \left( \theta \right) }{ \sqrt{3-2\,\cos \left( \theta \right) -2\,\sin \left( \theta \right) }}}\\
\displaystyle \cos \left( A \right) ={\frac {1-\sin \left( \theta \right) }{ \sqrt{3-2\,\cos \left( \theta \right) -2\,\sin \left( \theta \right) }}}\\
\arctan(z)+\arctan(1/z)=\pm \pi/2\,,
\end{align}
the sign of the latter chosen according to whether $\mathfrak{S}(z)>0$ or $\mathfrak{S}(z)<0$ respectively, where 
\begin{align} \nonumber
\mathfrak{S}(z) = & \mbox{  1} \mbox{ if } 0 < \Re(z) \mbox{ or } \Re(z) = 0 \mbox { and } 0 < \Im(z) \\
\mathfrak{S}(z) = & \mbox{ -1} \mbox{ if } 0 > \Re(z) \mbox{ or } \Re(z) = 0 \mbox { and } 0 > \Im(z).
\end{align}


The simplest of the special cases, corresponding to $s=0$, gives

\begin{equation}
\displaystyle \int_{0}^{\pi/2}\!{\frac {\sin \left( -\pi/2\,\sin \left( \theta \right) +\theta \right) {{\rm e}^{\pi/2\, \left( -1+\cos \left( \theta \right)  \right) }}}{\cosh \left( \pi\, \left( -1+\cos \left( \theta \right)  \right)  \right) 
\mbox{}+\cos \left( \pi\,\sin \left( \theta \right)  \right) }}+{\frac {\sin \left( \pi/2\,\sin \left( \theta \right) +\theta \right) {{\rm e}^{-\pi/2\, \left( -1+\cos \left( \theta \right)  \right) }}}{\cosh \left( \pi\, \left( -1+\cos \left( \theta \right)  \right)  \right) +\cos \left( \pi\,\sin \left( \theta \right)  \right) }}
\mbox{}\,{\rm d}\theta=1\,,
\label{As0}
\end{equation}
which, after the transformation $\theta\rightarrow\pi/2-\theta$ alternatively becomes
\begin{equation}
\displaystyle \int_{0}^{\pi/2}\!{\frac {\cos \left( \pi/2\,\cos \left( \theta \right) +\theta \right) {{\rm e}^{\pi/2\, \left( -1+\sin \left( \theta \right)  \right) }}+\cos \left( -\pi/2\,\cos \left( \theta \right) +\theta \right) {{\rm e}^{-\pi/2\, \left( -1+\sin \left( \theta \right)  \right) }}
\mbox{}}{\cosh \left( \pi\,(\sin \left( \theta \right) -1) \right) +\cos \left( \pi\,\cos \left( \theta \right)  \right) }}\,{\rm d}\theta=1\,.
\label{As0A}
\end{equation}
For comparison, see \eqref{Sa2cEp} with $n=0$, which, with $B_0=1$, can be rewritten
\begin{equation}
\displaystyle \int_{0}^{\pi/2}\!{\frac {\sin \left( \pi/2\,\cos \left( \theta \right) -\theta \right) {{\rm e}^{-\pi/2\,\sin \left( \theta \right) }}+\sin \left(\pi/2\,\cos \left( \theta \right) +\theta \right) {{\rm e}^{\pi/2\,\sin \left( \theta \right) }}
\mbox{}}{\cosh \left( \pi\,\sin \left( \theta \right)  \right) -\cos \left( \pi\,\cos \left( \theta \right)  \right) }}\,{\rm d}\theta=1
\label{Sa2Xa}
\end{equation}
or
\begin{equation}
\displaystyle \int_{0}^{\pi/2}\!{\frac {-\cos \left( \pi/2\,\sin \left( \theta \right) +\theta \right) {{\rm e}^{-\pi/2\,\cos \left( \theta \right) }}+\cos \left( -\pi/2\,\sin \left( \theta \right) +\theta \right) {{\rm e}^{\pi/2\,\cos \left( \theta \right) }}
\mbox{}}{\cosh \left( \pi\,\cos \left( \theta \right)  \right) -\cos \left( \pi\,\sin \left( \theta \right)  \right) }}\,{\rm d}\theta=1\,.
\label{Sa2Xb}
\end{equation}
Any of the above can also be rewritten in a different form after the transformation $\sin(\theta)=v$, giving for example the following variant of \eqref{Sa2cEp}:
\begin{equation}
\displaystyle \int_{0}^{1}\!{\frac {\sinh \left( \pi/2\, \sqrt{1-{v}^{2}} \right) \cos \left( \pi\,v/2 \right) }{\cosh \left( \pi\, \sqrt{1-{v}^{2}} \right) 
\mbox{}-\cos \left( \pi\,v \right) }}+{\frac {v\cosh \left( \pi/2\, \sqrt{1-{v}^{2}} \right) \sin \left( \pi\,v/2 \right) }{ \sqrt{1-{v}^{2}} \left( \cosh \left( \pi\, \sqrt{1-{v}^{2}} \right) -\cos \left( \pi\,v \right)  \right) 
\mbox{}}}\,{\rm d}v=1/2\,.
\label{Sa2XA}
\end{equation}

\subsection{Integration path C: Two straight lines: $v_{\pm}:=\pm i+(1\mp i)v $} \label{sec:twostraight}

By joining the integration endpoints $(\pm i,0)$ and $(0,1)$ with straight lines, (see Figure \ref{fig:Figure1}), \eqref{BasicEq} becomes

\begin{align} \nonumber
&\displaystyle {\frac {\zeta \left( s \right)  \left( 1-{2}^{1-s} \right) }{{2}^{s}}}=-\sin \left( s\pi/2 \right) \sum _{k=0}^{\infty }{E_{s}} \left(\pi\, \left( k+1/2 \right)  \right) \\ \nonumber
&+\int_{0}^{1}\! \left[ \frac {\cos \left( \pi v/2 \right) \cosh \left( \pi v/2 \right) +\sin \left( \pi v/2 \right) \sinh \left( \pi\,v/2 \right) 
}{\cos \left( \pi\,v \right) +\cosh \left( \pi\,v \right) }\cos ( s\arctan ( {\frac {v}{v-1}} )  ) \right.\\
&\left.-{\frac {\cos \left( \pi\,v/2 \right) \cosh \left( \pi v/2 \right) -\sin \left( \pi v/2 \right) \sinh \left( \pi v/2 \right) }{\cos \left( \pi\,v \right) +\cosh \left( \pi\,v \right) }\sin ( s\arctan ( {\frac {v}{v-1}} )  ) }
\mbox{} \right]  \left(  \left( v-1 \right) ^{2}+{v}^{2} \right) ^{-s/2}\,{\rm d}v\,.
\label{Linear}
\end{align}


The case $s=0$ yields the identity
\begin{equation}
\displaystyle \int_{0}^{1}\!{\frac {\cosh \left( \pi\,v/2\, \right) \cos \left( \pi\,v/2 \right) +\sinh \left( \pi\,v/2 \right) \sin \left( \pi\,v/2 \right) 
\mbox{}}{\cosh \left( \pi\,v \right) +\cos \left( \pi\,v \right) }}\,{\rm d}v=1/2\,,
\label{NeatId}
\end{equation}
which can also be rewritten in the interesting form

\begin{equation}
\displaystyle \int_{0}^{1}\!{\frac {\cosh \left( \pi\,v/2 \right) \cos \left( \pi\,v/2 \right) }{  \cosh^{2} \left( \pi\,v/2 \right) - \sin^{2} \left( \pi\,v/2 \right) 
\mbox{}}}\,{\rm d}v=1-\int_{0}^{1}\!{\frac {\sinh \left( \pi\,v/2 \right) \sin \left( \pi\,v/2 \right) }{  \cosh^{2} \left( \pi\,v/2 \right)  - \sin^{2} \left( \pi\,v/2 \right)}}\,{\rm d}t\,.
\label{NeatId2}
\end{equation}
{\bf Remarks:} After considerable execution time, the computer program Mathematica {\cite{Math}} successfully evaluates \eqref{NeatId}. Also, by unknown means, both Glasser (private communication) and Mathematica obtain the even more interesting variation:
\begin{equation}
\displaystyle \int_{0}^{1}\!{\frac {\cosh \left( \pi\,v/2 \right) \cos \left( \pi\,v/2 \right) }{  \cosh^{2} \left( \pi\,v/2 \right) - \sin^{2} \left( \pi\,v/2 \right) 
\mbox{}}}\,{\rm d}v=\frac{1}{2}+\frac{1}{\pi}\,\log\left(\coth(\pi/4)\right)\,.
\label{LG1}
\end{equation}

After subtracting \eqref{NeatId}, the case $s=-2$ yields
\begin{equation}
\displaystyle \int_{0}^{1}\!{\frac { \left( {v}^{2}-2\,v \right) \cos \left( \pi\,v/2 \right) \cosh \left( \pi\,v/2 \right) -\sinh \left( \pi\,v/2\right) \sin \left( \pi\,v/2 \right) {v}^{2}
}{\cosh \left( \pi\,v \right) +\cos \left( \pi\,v \right) }}\,{\rm d}v=-1/4\,.
\label{Sminus2}
\end{equation}

Similarly, the case $s=2$ produces the identity
\begin{equation}
\displaystyle \int_{0}^{1}\!{\frac {\cos \left( \pi\,v/2 \right) \cosh \left( \pi\,v/2 \right)  \left( 2\,{v}^{2}-1 \right) - \left( 2{v}^{2}-4\,v+1 \right) \sin \left( \pi\,v/2 \right) \sinh \left( \pi\,v/2 \right) 
\mbox{}}{ \left( \cosh \left( \pi\,v \right) +\cos \left( \pi\,v \right)  \right)  \left( 2\,{v}^{2}-2\,v+1 \right) ^{2}}}\,{\rm d}v=-{\pi}^{2}/48\,,
\label{Splus2}
\end{equation}

and, again after taking \eqref{NeatId} into consideration, the case $s=1$ gives
\begin{align}
\displaystyle \int_{0}^{1}\;&\frac {v \left( v-1 \right) \cos \left( \pi\,v/2 \right) \cosh \left( \pi\,v/2 \right) +\sinh \left( \pi\,v/2 \right) \sin \left( \pi\,v/2 \right) {v}^{2}
\mbox{}}{ \left( \cosh \left( \pi\,v \right) +\cos \left( \pi\,v \right)  \right)  \left( 2\,{v}^{2}-2\,v+1 \right) }\,{\rm d}v=\\
&-1/2\,\sum _{k=0}^{\infty }{E_1} \left(\pi\,(k+1/2) \right) -\ln  \left( 2 \right)/4 
\mbox{}+1/4
\label{seq1}
\end{align}
Note that in all of the above, the denominator factor $\left( 2\,{v}^{2}-2\,v+1 \right)\neq 0$ on the real line.\newline

\subsection{Integration Path D: Four straight lines}


In this case, the integration endpoints of \eqref{BasicEq} are joined by straight lines (see Figure \ref{fig:Figure1}) in the following order: $(0,-i)$, $(1,-i)$, $(1,0)$, $(1,+i)$ and $(0,+i) $ 
and, after the appropriate change of variables and simplification, the following identity emerges:
\begin{align} \nonumber
&\displaystyle \zeta \left( s \right)  \left( 1-{2}^{1-s} \right) =-{2}^{s}\sin \left( s\pi/2 \right) \sum _{k=0}^{\infty }{E_{s}} \left(\pi\, \left( k+1/2 \right)  \right) \\ \nonumber
&+{2}^{s-1} \int_{0}^{1}\! \left({\frac {\cos \left( -\pi\,v/2+s\arctan \left( 1/{v} \right)  \right) {{\rm e}^{\pi/2}}-\cos \left( \pi\,v/2+s\arctan \left( 1/{v} \right)  \right) {{\rm e}^{-\pi/2}}
\mbox{}}{\cosh \left( \pi \right)-\cos \left( \pi\,v \right)  }}\right.  \\
& \left.  \hspace{60pt} +  {\frac {\cos \left( s\arctan \left( v \right)  \right) }{\cosh \left( \pi\,v/2 \right) }} \right)  \left( {v}^{2}+1 \right) ^{-s/2}\,{\rm d}v \,.
\label{Path4}
\end{align}
As before, many special cases are embedded in \eqref{Path4}. In the case $s=0$, the simplified form involves integrals all of which are amenable to analysis by both Mathematica and Maple \cite{Math}, \cite{Maple} and so this case is not discussed here. Similar but simpler integrals can be found in \cite[Sections 3.531 and 3.532]{G&R}.

\section{Summary}
What is believed to be a new series representation of Riemann's zeta function has been developed, from which many novel identities flowed, most of which appear to have no obvious application. It is suggested that further consideration of the integrals $J_{1}(t)$ and $J_{2}(t)$ in \eqref{ZmRSoln} and \eqref{ZmISoln} may shed light on the behaviour of $\zeta^{-}(1/2+it)$ in the asymptotic limit $t\rightarrow \infty$.

\section{Acknowledgements}
The author thanks Larry Glasser for making him aware of Russell's 1876 work \cite{RussellA} and Christophe Vignat for pointing out \eqref{EuInt}.

\bibliographystyle{unsrt}

\bibliography{biblio}
\begin{appendices}

\section{Notation} \label{sec:AppendixA}

Throughout, the symbols $j,k,m,n,p,r,N$ are positive integers, non-zero except as noted. $E_s(z)$ is the (generalized) exponential integral, $\mathcal{E}(n)$ is an Euler number, $\mathcal{E}(n,z)$ is an Euler polynomial, $B_n$ is a Bernoulli number, $H_n$ is a harmonic number, $erf$ is the error function, $\gamma$ is the Euler-Mascheroni constant and $Si(z)$ is the sine integral \eqref{SiDef}. The floor function (``greatest integer less than") is symbolized as usual by $\lfloor...\rfloor$. Subscripted $\zeta_{R}(s)$ and $\zeta_{I}(s)$ refer to the real and imaginary components of $\zeta(s)$ or its incomplete forms. Other symbols are defined where they first appear. Any summation where the lower summation limit exceeds the upper summation limit is identically zero. At times, a complicated computer simplification has been used \cite{Maple} to reduce an expression; such use is indicated by the notation ``(Maple)".

\section{$\zeta(s)$ represented by various sums involving $E_s(z)$} \label{sec:EiSums}

In the literature, a number of representations are known in which $\zeta(s)$ or its associate $\xi(s)$ are expressed in terms of exponential integral functions. 
\begin{itemize}
\item
From Paris \cite{ParisExp} (where $\lambda$ is a parameter satisfying $|\arg(\lambda)|\leq \pi/2$) we have, with modified notation,
\begin{align}
\displaystyle \zeta \left( s \right) ={\frac {{(\pi\,\lambda)}^{s/2}}{\Gamma \left( s/2 \right) } \left( \sum _{n=1}^{\infty }{E_{1-s/2}} \left(\pi\,{n}^{2}\lambda \right) +{\frac {1}{ \sqrt{\lambda}}\sum _{n=1}^{\infty }{E_{s/2+1/2}} \left({\frac {\pi\,{n}^{2}}{\lambda}} \right) }
\mbox{}+{\frac {1}{ \sqrt{\lambda} \left( s-1 \right) }}-\frac{1}{s} \right) }\,.
\label{Paris}
\end{align}
\item
From \cite{LeClair:2013oda} we have the representation
\begin{align}
\xi(s)= \pi \left( s-1 \right) \sum _{n=1}^{\infty }{n}^{2}{E_{-s/2}} \left(\pi\,{n}^{2} \right)-\pi\,s\sum _{n=1}^{\infty }{n}^{2}{ E_{(s-1)/2}} \left( \pi\,{n}^{2} \right)+4\,\pi\,\sum _{n=1}^{\infty }{n}^{2}{{\rm e}^{-\pi\,{n}^{2}}}\,,
\mbox{} 
\label{Leclair1}
\end{align}
where
\begin{equation}
\xi(s)\equiv \displaystyle { \left( s-1 \right) {\pi}^{-s/2}}{}\,\zeta \left( s \right) \Gamma \left(1+ s/2 \right)
\label{xidef}
\end{equation}

and where it was later shown (\cite[Eq.(3)]{Romik}) that 
\begin{equation}
\displaystyle 4\pi\sum _{n=1}^{\infty }{n}^{2}{{\rm e}^{-\pi\,{n}^{2}}}={\frac {\pi^{1/4}}{2\Gamma \left( 3/4 \right) }}\,.
\label{Larry1}
\end{equation}
\end{itemize}
\section{The incomplete functions $\zeta^{\pm}(s)$ and $\eta^{\pm}(s)$, definitions and Lemmas } \label{sec:Def_and_Lemmas}

In  this Appendix, various results needed in the text are either defined, proven, or collected from the literature.
\begin{itemize}
\item
The sum $\displaystyle \sum _{k=0}^{\infty }{\it E_s} \left(\pi\, \left( k+1/2 \right)  \right)$ appears frequently. Some identities pertaining to this series follow. From the integral representation \eqref{ExpInt}, interchange the series and integration to obtain
\begin{equation}
\displaystyle \sum _{k=0}^{\infty }{\it E_s} \left(\pi\, \left( k+1/2 \right)  \right) =\frac{1}{2}\,\int_{1}^{\infty }\!{\frac {{v}^{-s}}{\sinh \left( \pi\,v/2 \right) }}\,{\rm d}v\,.
\label{Esum}
\end{equation}

From \cite[Eq.(29)]{Milgram} we have the representation
\begin{equation}
\displaystyle \zeta \left( s \right) ={\frac {{2}^{s-1}\,\sin \left(\pi\,s/2 \right) 
\mbox{}}{{2}^{1-s}-1}\int_{0}^{\infty }\!{\frac {{v}^{-s}}{\sinh \left( \pi\,v /2\right) }}\,{\rm d}v} \hspace{12pt} \sigma<0,
\label{ZetaEq29}
\end{equation}
after a simple change of variables. Define the upper and lower incomplete functions $\zeta^{\pm}(s)$ as

\begin{equation}
\displaystyle \zeta^{+} \left( s \right) ={\frac {{2}^{s-1}\,\sin \left(\pi\,s/2 \right) 
\mbox{}}{{2}^{1-s}-1}\int_{1}^{\infty }\!{\frac {{v}^{-s}}{\sinh \left( \pi\,v /2\right) }}\,{\rm d}v} \,,
\label{ZetaPlus}
\end{equation}
and
\begin{equation}
\displaystyle \zeta^{-} \left( s \right) ={\frac {{2}^{s-1}\,\sin \left(\pi\,s/2 \right) 
\mbox{}}{{2}^{1-s}-1}\int_{0}^{1 }\!{\frac {{v}^{-s}}{\sinh \left( \pi\,v /2\right) }}\,{\rm d}v} \hspace{12pt} \sigma<0,
\label{ZetaMinus}
\end{equation}

and similarly for the alternating function
\begin{equation}
\eta^{\mp}(s)\equiv (1-2^{1-s})\zeta^{\mp}(s)\,
\label{etaPmDef}
\end{equation}
so that
\begin{equation}
\zeta(s)=\zeta^{-}(s)+\zeta^{+}(s)
\label{zAdd}
\end{equation}
and similarly for $\eta(s)$. Thus, \eqref{Esum} becomes
\begin{equation}
\displaystyle \sum _{k=0}^{\infty }{\it E_s} \left(\pi\, \left( k+1/2 \right)  \right) =\frac{1}{2}\int_{1}^{\infty }\!{\frac {{v}^{-s}}{\sinh \left( \pi\,v/2 \right) }}\,{\rm d}v={\frac { \zeta^{+}\left( s \right)  \left( {2}^{1-s}-1 \right) }{{2}^{s}\,\sin \left( \pi\,s/2 \right) 
}}=-{\frac { \eta^{+}\left( s \right)   }{{2}^{s}\,\sin \left( \pi\,s/2 \right) 
}}.
\label{Zplus}
\end{equation}

{\bf Remark:} By writing \eqref{ZetaMinus} as follows


\begin{equation}
\displaystyle \zeta^{-}(s)= \,{\frac {{2}^{s-1}\sin \left( \pi\,s/2 \right) }{{2}^{1-s}-1
\mbox{}} \left( \int_{0}^{1}\!\left( {\frac {{v}^{-s}}{\sinh \left( \pi\,v/2 \right) }}-{\frac {2\,{v}^{-s-1}}{\pi}}\,{\rm d}v\right)-{\frac {2}{\pi\,s}} \right) }\,\hspace{10pt}\sigma<1\,,
\label{redefine}
\end{equation}
we find
\begin{align} \label{Zm0}
\zeta^{-}(0) &=-1/2  \\
\zeta^{+}(0) & =\hspace{12pt}0\,.
\label{Zp0}
\end{align}

Similarly, by expanding \eqref{ZetaEq29} about $s=1$, from the (leading) term of order $(s-1)^{-1}$ we identify 

\begin{equation}
\displaystyle \int_{0}^{\infty }\!\left({\frac {1}{v\sinh \left( \pi\,v/2 \right) }}-{\frac {2}{\pi\,{v}^{2}}}\right)\,{\rm d}v=-\ln  \left( 2 \right)
\label{Zs1}
\end{equation}

and from the term of order $(s-1)^{0}$ with known coefficient $\gamma$ we find

\begin{equation}
\displaystyle \int_{0}^{\infty }\!{\frac {\ln  \left( v \right) }{v} \left(  \frac{1}{ \sinh \left( \pi\,v /2\right)  } -{\frac {2}{\pi\,v}} \right) }\,{\rm d}v= \left( \gamma-3/2\,\ln  \left( 2 \right)  \right)
\mbox{}\ln  \left( 2 \right)\,.
\label{Zs2}
\end{equation}

\item
From Erdelyi \cite[Eq. 1.5.1(37)]{Erdelyi1} with $-1<-\Re(s)-2n<0$:
\begin{equation}
\displaystyle \int_{0}^{\infty }\!{v}^{-s-2\,n}\cos \left( \pi\, \left( k+1/2 \right) v \right) \,{\rm d}v= \left( \pi\, \left( k+1/2 \right)  \right) ^{s+2\,n-1}
\mbox{}\Gamma \left( -s-2\,n+1 \right)  \left( -1 \right) ^{n}\sin \left( \pi\,s/2 \right)\,. 
\label{Erd37}
\end{equation}

\item
From \cite[Eq. (2.21) with $j=0$]{Milgram:1985}:
\begin{equation}
\displaystyle { E_{-p}} \left(\pm i\pi\, \left( j+1/2 \right)  \right) ={\frac {\Gamma \left( p+1 \right) {{\rm e}^{\mp i\pi\, \left( j+1/2 \right) }}}{ \left(\pm i\,\pi\, \left( j+1/2 \right)  \right) ^{p+1}}\sum _{k=0}^{p}{\frac { \left( \pm \,i\pi\, \left( j+1/2 \right)  \right) ^{k}}{\Gamma \left( k+1 \right) }}}\,.
\label{EiNegP}
\end{equation}

\item
From Hansen \cite[Eqs. (17.4.10) and (14.4.10) respectively]{Hansen} or NIST \cite[Eqs. (24.8.4) and (24.8.5)]{NIST} with $0\le v/2\le 1$:
\begin{equation}
\displaystyle \sum _{k=0}^{\infty }{\frac {\cos \left( \pi\, \left( k+1/2 \right) v \right) }{ \left( k+1/2 \right) ^{2\,n}}}={\frac { \left( 2\,\pi \right) ^{2\,n} \left( -1 \right) ^{n}
\mbox{}{\bf { \mathcal{E}}} \left( 2\,n-1,v/2 \right) }{4\,\Gamma \left( 2\,n \right) }}
\label{hans1}
\end{equation}

and

\begin{equation}
\displaystyle \sum _{k=0}^{\infty }{\frac {\sin \left( \pi\, \left( k+1/2 \right) v \right) }{ \left( k+1/2 \right) ^{2\,n+1}}}={\frac { \left( 2\,\pi \right) ^{2\,n+1}
\mbox{} \left( -1 \right) ^{n}\mathcal{E} \left( 2\,n,v/2 \right) }{4\Gamma \left( 2\,n+1 \right) }}
\label{Hans2}
\end{equation}
where $\mathcal{E} \left( 2\,n-1,v/2 \right)$ represents the Euler polynomial, defined (from \cite[FunctionAdvisor]{Maple}) in terms of Bernoulli numbers $B_k$ by 

\begin{equation}
\displaystyle \mathcal{E} \left( n,z \right) =\Gamma \left( n+1 \right) \sum _{k=0}^{n}{\frac {{z}^{k}{B_{1-k+n}}  }{\Gamma \left( k+1 \right) \Gamma \left( 2-k+n \right) } \left( 2- {2}^{2-k+n}  \right) }
\label{EulerNo}
\end{equation}
which, reversed, becomes
\begin{equation}
\displaystyle \mathcal{E} \left( n,z \right) =\Gamma \left( n+1 \right) \sum _{k=0}^{n}{\frac { \left( 2-{2}^{2+k} \right) {z}^{n-k}{B_{k+1}} }{\Gamma \left( 1-k+n \right) \Gamma \left( 2+k \right) 
\mbox{}}}\,,
\label{EulernoRev}
\end{equation}

reducing to
\begin{equation}
\displaystyle {\mathcal{E}} \left( n \right) =\Gamma \left( n+1 \right) \sum _{k=0}^{n}{\frac { \left( {2}^{k+1}-{2}^{2+2\,k} \right) B_{k+1} }{\Gamma \left( 1-k+n \right) 
\mbox{}\Gamma \left( k+2 \right) }}\,,
\label{EuNo}
\end{equation}
if $z=1/2$ because of the definition \cite[Eq. 24.2.9]{NIST}
\begin{equation}
\mathcal{E}(n)=2^n\mathcal{E}(n,1/2)\,.
\label{Edef}
\end{equation}
The following useful identity \cite[Eq. (24.4.15]{NIST} inverts the above: 
\begin{equation}
\displaystyle B_{2n} ={\frac {\Gamma \left( 2\,n+1 \right) }{{2}^{2\,n} \left( {2}^{2\,n}-1 \right) }\sum _{k=0}^{n-1}{\frac {\mathcal{E} \left( 2\,k \right) }{\Gamma \left( 2\,k+1 \right) \Gamma \left( 2\,n-2\,k \right) }}}\,.
\label{B2n}
\end{equation}

\item
From \cite[Eq. 24.2.10]{NIST}

\begin{equation}
\displaystyle \mathcal{E} \left( j,\frac{1}{2x}+\frac{1}{2} \right) =\sum _{k=0}^{j}\binom{j}{ k}\frac {\mathcal{E} \left( k \right) }{{2}^{k}} \left(\frac{1}{2x} \right) ^{j-k}\,,
\label{EuIdhalf}
\end{equation}
we have the following lemma\newline

{\bf Lemma:} For $m>0$ and $t\geq 1$,

\begin{equation}
\displaystyle \Gamma \left( 2\,m+2 \right) \sum _{j=0}^{2\,m} \sum _{k=0}^{j}{\frac { \left( -1 \right) ^{j}\mathcal{E} \left( k \right) {t}^{-k}}{\Gamma \left( k+1 \right) \Gamma \left( j-k+1 \right) 
\mbox{} \left( j+1 \right) \Gamma \left( 2\,m-j+1 \right) }}  =\sum _{k=0}^{2\,m} \left( -1 \right) ^{k}{t}^{-k}\mathcal{E} \left( k \right) \,.
\label{Lemma1}
\end{equation}
\newline

{\bf Proof:}

Reorder the series in \eqref{Lemma1} giving

\begin{align} \nonumber
 \displaystyle \sum _{j=0}^{2\,m}  &\sum _{k=0}^{j}{\frac { \left( -1 \right) ^{j}\mathcal{E} \left( k \right) {t}^{-k}}{\Gamma \left( k+1 \right) \Gamma \left( j-k+1 \right) 
\mbox{} \left( j+1 \right) \Gamma \left( 2\,m-j+1 \right) }} \\
&  =\sum _{k=0}^{2\,m}  {\frac { \left( -1 \right) ^{k}
\mbox{}{t}^{-k}\mathcal{E} \left( k \right) }{\Gamma \left( k+1 \right) }\sum _{j=0}^{2\,m-k}{\frac { \left( -1 \right) ^{j}}{\Gamma \left( j+1 \right)  \left( j+k+1 \right) \Gamma \left( 2\,m-j-k+1 \right) }}}\,. 
\label{LemProof1}
\end{align}

The inner sum can be identified

\begin{equation}
\displaystyle \sum _{j=0}^{2\,m-k}{\frac { \left( -1 \right) ^{j}}{\Gamma \left( j+1 \right)  \left( j+k+1 \right) \Gamma \left( 2\,m-j-k+1 \right) }}={\frac {{\mbox{$_2$F$_1$}(k+1,-2\,m+k;\,k+2;\,1)}}{ \left( k+1 \right) 
\mbox{}\Gamma \left( 2\,m+1-k \right) }}
\label{LemProof2}
\end{equation}
and the hypergeometric function evaluates to 
\begin{equation}
\displaystyle {\mbox{$_2$F$_1$}(k+1,-2\,m+k;\,k+2;\,1)}={\frac {\Gamma \left( k+2 \right) \Gamma \left( 2\,m+1-k \right) }{\Gamma \left( 2\,m+2 \right) }}\,,
\label{LemProof3}
\end{equation}

all of which is easily seen to reduce to \eqref{Lemma1}.  {\bf Q.E.D.}\newline

\item
From \cite{DirBeta}, a sum appearing frequently in the text, is identified by any of the following: 
\begin{equation}
\displaystyle   \sum _{k=0}^{\infty }{\frac { \left( -1 \right) ^{k}}{ \left( k+1/2 \right) ^{2\,n+1}}}={\frac {\psi \left( 2\,n,3/4 \right) -\psi \left( 2\,n,1/4 \right) }{\Gamma \left( 2\,n+1 \right) 
{2}^{2\,n+1}}}  =\beta \left( 2\,n+1 \right) {2}^{2\,n+1}={\frac { \left( -1 \right) ^{n}{\pi}^{2\,n+1}
\mathcal{E} \left( 2\,n \right) }{2\Gamma \left( 2\,n+1 \right) }}
\label{BetaDef}
\end{equation}
where $\beta(x)$ is the Dirichlet Beta Function and $\psi(n,x)$ is the polygamma function.\newline

\item
From Mathematica \cite{Math} we have the following sum

\begin{align} \nonumber
\displaystyle \sum _{j=0}^{\infty }{\frac {\psi \left( 2\,j+1 \right)  \left( -1 \right) ^{j}{z}^{-j}{2}^{2\,j}}{\Gamma \left( 2\,j+1 \right) }}& = -\left( { \rm Ci} \left( \frac{2}{\sqrt{z}} \right) +\frac{\ln  \left( z \right)}{2}  \right) \cos \left(  \frac{2}{\sqrt z}  \right) 
\mbox{}\\&
+\ln  \left( 2 \right) \cosh \left(\frac{2}{\sqrt{-z}}  \right)
-\sin \left( \frac{2}{\sqrt{z}}  \right) {\rm Si} \left( \frac{2}{\sqrt{z}} \right) 
\label{Math1a}
\end{align}
where $ \rm Ci$ and $ \rm Si$ are the cosine and sine integrals respectively (see \cite[Eq. 6.2.9 and 6.2.16]{NIST}). In the case that $ z=-4/\pi^2/(k+1/2)^2 $, \eqref{Math1a} reduces to
\begin{equation}
\displaystyle \sum _{j=0}^{\infty }{\frac {\psi \left( 2\,j+1 \right)  \left( -1 \right) ^{j} \left( k+1/2 \right) ^{2\,j}{\pi}^{2\,j}
\mbox{}}{\Gamma \left( 2\,j+1 \right) }}=- \left( -1 \right) ^{k}{\rm Si} \left( \pi(k+1/2) \right) \,.
\label{Math1b}
\end{equation}

\item
Taking \eqref{EiNeg_m} into account, the series appearing in \eqref{Sa2cOm} with $n>0$, can be written


\begin{equation}
\displaystyle \sum _{k=0}^{\infty }{E_{1-2n}} \left( \pi\, \left( k+1/2 \right)  \right) ={\frac {\Gamma \left( 2\,n \right) }{\pi}\sum _{j=0}^{2\,n-1}  {\frac {{{\pi^{-j}\,\rm e}^{-\pi/2}}}{\Gamma \left( 2\,n-j \right)}\sum _{k=0}^{\infty }{\frac {{{\rm e}^{-\pi\,k}}}{ \left( k+1/2
\mbox{} \right) ^{j+1}}}}  }\,,
\label{EiSumNeg}
\end{equation}
and from the standard definition \cite{Lerch} of the Lerchphi transcendent $\Phi$, we have
\begin{equation}
\displaystyle \sum _{k=0}^{\infty }{\frac {{{\rm e}^{-\pi\,k}}}{ \left( k+1/2 \right) ^{j+1}}}={\Phi} \left( {{\rm e}^{-\pi}},j+1,1/2 \right) \,
\label{LerchPhi}
\end{equation}
which can be rewritten in terms of generic polylog functions $Li_{j}(x)$ using the identities

\begin{equation}
\displaystyle {\Phi} \left( z,s,a \right) ={\frac {{\Phi} \left( {z}^{2},s,a/2 \right) +z{\Phi} \left( {z}^{2},s,a/2+1/2 \right) }{{2}^{s}}}
\label{LerchId}
\end{equation}
with $a=1$ and \cite[Eq. (6)]{Lerch}, eventually, with $n>0$, yielding the identity 

\begin{equation}
\displaystyle \sum _{k=0}^{\infty }{E_{1-2n}} \left(\pi\, \left( k+1/2 \right)  \right) ={\frac {\Gamma \left( 2\,n \right) }{\pi} \left( \sum _{j=0}^{2\,n-1}{\frac { Li_{j+1} \left(  {{\rm e}^{-\pi/2}} \right)   {2}^{j+1}}{\Gamma \left( 2\,n-j \right) 
\mbox{}{\pi}^{j}}}-\sum _{j=0}^{2\,n-1}{\frac {Li_{j+1} \left( {{\rm e}^{-\pi}}  \right) }{\Gamma \left( 2\,n-j \right) {\pi}^{j}}}
\mbox{} \right) }\,.
\label{EiPoly}
\end{equation}
{\bf Remark:} Relevant to the above, note both the identity
\begin{equation}
Li_{1}(x)=-ln(1-x)
\label{LiLn}
\end{equation}
and the connection to the Bose-Einstein distribution (see \cite[Eq. (13)]{Lerch}).

\end{itemize}

\section{Appendix C - a Theorem} \label{sec:Theorem}

In the main text, we are interested in combinations of the exponential integral function of complex conjugate argument. Here we present a proof of the following identity, where $k\in 0,1\dots$, $p\in 1,2\dots$ and $\rm{Si}(z)$ is the sine integral:
\begin{equation}
\rm{Si}(z)= {\int_{0}^{\, z} \frac{\sin(t)}{t}\, \rm{d}t}\,.
\label{SiDef}
\end{equation}

\begin{theorem}
{\rm For integer $p\geq 1$,}
\begin{align} \nonumber
\displaystyle \left( -1 \right) ^{p}\Gamma \left( 2\,p \right)& \frac {    { E_{2p}} \left( -i\pi\, \left( k+1/2 \right)  \right) + {E_{2p}} \left( i\pi\, \left( k+1/2 \right)  \right)  
\mbox{}} { \left( \pi\, \left( k+1/2 \right)  \right) ^{2\,p-1}}\\
 &=2\, \left( -1 \right) ^{k}\sum _{j=1}^{p-1}{\frac { \left( -1 \right) ^{j}\Gamma \left( 2\,j \right) }{ \left( \pi\, \left( k+1/2 \right)  \right) ^{2\,j}}} 
- 2 \,(\rm{Si} (\pi  (k + 1/2))-\pi/2 )
\label{Gans4}
\end{align}
\end{theorem}

{\bf Proof:}\newline

Begin by converting each of the functions $E_{2p}$ and $\rm Si$ to a (convergent) integral representation using \eqref{ExpInt} and \eqref{SiDef} respectively, giving an equivalent form of \eqref{Gans4}

\begin{align} \nonumber
&\displaystyle \int_{1}^{\infty }\!{v}^{-2\,p}\cos \left(  \left( k+1/2 \right) v\pi \right) \,{\rm d}v\\=
&{\frac { \left( -1 \right) ^{p} \left( \pi\, \left( k+1/2 \right)  \right) ^{2\,p-1}}{\Gamma \left( 2\,p \right) } \left( -\int_{0}^{1}\!{\frac {\sin \left(  \left( k+1/2 \right) v\pi \right) }{v}}\,{\rm d}v
\mbox{}+ \left( -1 \right) ^{k}\sum _{j=1}^{p-1}{\frac { \left( -1 \right) ^{j}\Gamma \left( 2\,j \right) }{ \left( k+1/2 \right) ^{2\,j}{\pi}^{2\,j}}}
\mbox{}+\pi/2 \right) }
\label{Pit}
\end{align}

After integrating once by parts, the left-hand side of \eqref{Pit} becomes

\begin{equation}
\displaystyle \int_{1}^{\infty }\!{v}^{-2\,p}\cos \left(  \left( k+1/2 \right) v\pi \right) \,{\rm d}v=-{\frac { \left( -1 \right) ^{k}}{\pi\, \left( k+1/2 \right) }}
\mbox{}+{\frac {2p}{\pi\, \left( k+1/2 \right)
\mbox{}}\int_{1}^{\infty }\!\sin \left(  \left( k+1/2 \right) v\pi \right) {v}^{-1-2\,p}\,{\rm d}v}\,.
\label{PitEq}
\end{equation}

We proceed by induction. Set $p=1$ in \eqref{Pit}, integrate again by parts, complete the incomplete integral by utilizing
\begin{equation}
\int_{0}^{\infty}\frac{\sin(t)}{t}\rm{d}t = \frac{\pi}{2}
\end{equation}

 and simplify (Maple) to find that \eqref{Pit} becomes
 
\begin{align} \nonumber
&\displaystyle 2\,\pi\, \left( k+1/2 \right) {\rm Si} \left( \pi\, \left( k+1/2 \right)  \right) -{\frac { 2\,\left( -1 \right) ^{k}}{\pi\, \left( k+1/2 \right) }}+{\frac {2\, \left( -1 \right) ^{k}-{\pi}^{3} \left( k+1/2 \right) ^{2}}{\pi\, \left( k+1/2 \right) }}
\mbox{}\\
=&2\,\pi\, \left( k+1/2 \right) {\rm Si} \left( \pi\, \left( k+1/2 \right)  \right) -{\pi}^{2} \left( k+1/2 \right) 
\label{CasePeq1}
 \end{align}
 which is obviously true. Having shown that it is true for $p=1$, use \eqref{PitEq} to posit that the equivalent general form of \eqref{Gans4} is true for any value of $p$:
\begin{align} \nonumber
&\displaystyle \int_{1}^{\infty }\!\sin \left(  \left( k+1/2 \right) v\pi \right) {v}^{-1-2\,p}\,{\rm d}v
\mbox{}=-{\frac { \left( -1 \right) ^{p} \left( \pi\, \left( k+1/2 \right)  \right) ^{2\,p}}{\Gamma \left( 2\,p+1 \right) }\int_{0}^{1}\!{\frac {\sin \left(  \left( k+1/2 \right) v\pi \right) }{v}}\,{\rm d}v}\\
&+{\frac { \left( -1 \right) ^{k+p} \left( \pi\, \left( k+1/2 \right)  \right) ^{2\,p}}{\Gamma \left( 2\,p
\mbox{}+1 \right) }\,\sum _{j=1}^{p-1}{\frac { \left( -1 \right) ^{j}\Gamma \left( 2\,j \right) }{ \left( \pi\, \left( k+1/2 \right)  \right) ^{2\,j}}}}
\mbox{}+\left( -1 \right) ^{p}\frac{\pi}{2}{\frac { \left( \pi\, \left( k+1/2 \right)  \right) ^{2\,p}}{\Gamma \left( 2\,p+1 \right) }}+{\frac { \left( -1 \right) ^{k}}{2\,p}}\,;
\label{posit}
 \end{align}

then transform \eqref{posit} by setting $p\rightarrow p+1$. Integrating by parts converts the left-hand side of the transformed \eqref{posit} as follows:
\begin{equation}
\displaystyle \int_{1}^{\infty }\!\sin \left(  \left( k+1/2 \right) v\pi \right) {v}^{-3-2\,p}\,{\rm d}v
\mbox{}=-{\frac { \left( k+1/2 \right) ^{2}{\pi}^{2}\int_{1}^{\infty }\!\sin \left(  \left( k+1/2 \right) v\pi \right) {v}^{-1-2\,p}\,{\rm d}v
\mbox{}}{ 2\left( 2\,p+1 \right)  \left( 1+p \right) }}+{\frac { \left( -1 \right) ^{k}}{2+2\,p}}\,.
\label{Sintp1}
\end{equation}
Similarly, convert the transformed sum within \eqref{posit} as follows
\begin{equation}
\displaystyle \sum _{j=1}^{p}{\frac { \left( -1 \right) ^{j}\Gamma \left( 2\,j \right) }{ \left( k+1/2 \right) ^{2\,j}{\pi}^{2\,j}}}=\sum _{j=1}^{p-1}{\frac { \left( -1 \right) ^{j}\Gamma \left( 2\,j \right) }{ \left( k+1/2 \right) ^{2\,j}{\pi}^{2\,j}}}
\mbox{}+{\frac { \left( -1 \right) ^{p}\Gamma \left( 2\,p \right) }{ \left( k+1/2 \right) ^{2\,p}{\pi}^{2\,p}}}\,,
\label{NewSum}
\end{equation}
substitute both \eqref{Sintp1} and \eqref{NewSum} into the transformed \eqref{posit}, simplify (Maple) and it will be seen that, under the change $p\rightarrow p+1$, the transformed \eqref{posit} reduces to itself. By the logic of induction, \eqref{posit} is true and therefore so is \eqref{Gans4}. {\bf QED} \newline

{\bf Remark:} Since integration-by-parts underlies the recursion rule \eqref{recur}, the proof of \eqref{Gans4} could have been alternatively obtained by the $p-$fold application of \eqref{recur}.\newline

Since each of the integrals in \eqref{posit} can be expressed as a hypergeometric function (see \cite[Eqs. 6.2(5) and 6.2(6)]{Luke}) we find\newline

{\bf Corollary}

\begin{align} \nonumber
 \mbox{$_1$F$_2$}(\frac{1}{2}-p\,;\,\frac{1}{2},\frac{3}{2}-p\,;&\,-\frac{{\pi}^{2}\, \left( k+1/2 \right) ^{2}}{4})={\frac { \left( -1 \right) ^{p+1} \left( \pi\, \left( k+1/2 \right)  \right) ^{2\,p}
\mbox{}}{\Gamma \left( 2\,p-1 \right) }}{\mbox{$_1$F$_2$}(\frac{1}{2};\,\frac{3}{2},\frac{3}{2};\frac{-\pi^{2}\,\left( k+1/2 \right)^{2}}{4} )}\\
&+{\frac { \left( -1 \right) ^{k+p} \left( \pi\, \left( k+1/2 \right)  \right) ^{2\,p-1}}{\Gamma \left( 2\,p
\mbox{}-1 \right) }\sum _{j=1}^{p-1}{\frac {\Gamma \left( 2\,j \right)  \left( -1 \right) ^{j}}{ \left( \pi\, \left( k+1/2 \right)  \right) ^{2\,j}}}}\,.
\label{Corollary1}
\end{align} 

\section{A lengthy derivation} \label{sec:Circuitous}

This Appendix contains the details of a lengthy analysis based on \eqref{Cp}. With reference to that identity, and to reiterate: both sides are independent of the integer $p>0$. Consider the final term in \eqref{Cp}:

\begin{equation} 
- \left( -1 \right) ^{m}{\pi}^{2\,m}\sum _{j=0}^{p-2}{\frac {\Gamma \left( 2\,j+2 \right) \mathcal{E} \left( 2\,j+2+2\,m \right) }{\Gamma \left( 2\,j+3+2\,m \right) }} \,.
\label{Last}
\end{equation}
For large values of $n$ (and therefore $p$), an asymptotic approximation to the Euler numbers $\mathcal{E}(n)$ can be extracted by setting  $z=1/2$ and truncating (at $N$ terms) the following expression \cite[Eq. (3.6)]{LopezTemme} for the Euler polynomial $\mathcal{E}(n,z)$
\begin{equation}
\displaystyle \mathcal{E}(n,z)=4\,\Gamma \left( n+1 \right) \sum _{k=0}^{\infty }{\frac {\sin \left(  \left( 2\,k+1 \right) \pi\,z-\pi\,n/2 \right) }{ \left(  \left( 2\,k+1 \right) \pi \right) ^{n+1}}}\,,
\label{LopzTem}
\end{equation}
along with the definition \eqref{Edef} giving  

\begin{equation}
\displaystyle \mathcal{E} \left( 2\,j+2+2\,m \right) \approx -{\frac {2\,\Gamma \left( 2\,j+3+2\,m \right) 
\mbox{} \left( -1 \right) ^{j+m}}{{\pi}^{2\,j+3+2\,m}}\sum _{k=0}^{N}{\frac { \left( -1 \right) ^{k}}{\quad \left( k+1/2 \right) ^{2\,j+3+2\,m}}}}
\label{Eu2Id}
\end{equation}
where $N=1,2,\dots$. In the case $ \underset{N=\infty}\lim$, \eqref{Eu2Id} is equivalent to \eqref{BetaDef} and the approximation becomes an equality. In general, and for later use in the derivation of \eqref{Cpa4}, note that

\begin{align} \nonumber
\displaystyle 2\,\sum _{j=1}^{p-1}   \left( -1 \right) ^{p-j}\Gamma \left( 2\,p-2\,j \right) {\pi}^{2\,j}&\sum _{k=0}^{\infty } \frac{\left( -1 \right) ^{k}}{ \left( k+1/2 \right) ^{2\,p+2\,m-2\,j+1}}
\mbox{} \\
&=\pi\, \left( -1 \right) ^{m}{\pi}^{2\,p+2\,m}\sum _{j=0}^{p-2}{\frac {\Gamma \left( 2\,j+2 \right) \mathcal{E} \left( 2\,j+2+2\,m \right) }{\Gamma \left( 2\,j+3+2\,m \right) }}\,.
\label{Fsum1a}
\end{align}

For the moment however, we shall treat \eqref{Eu2Id} as an approximation - in fact it is a very good one numerically even for small values of $N$. When \eqref{Eu2Id} is substituted into \eqref{Last}, both sums are finite and can therefore be interchanged, yielding

\begin{equation}
\displaystyle \left( -1 \right) ^{m}{\pi}^{2\,m}\sum _{j=0}^{p-2}{\frac {\Gamma \left( 2\,j+2 \right) \mathcal{E} \left( 2\,j+2+2\,m \right) }{\Gamma \left( 2\,j+3+2\,m \right) }}
\mbox{}\approx-\,{\frac {2}{{\pi}^{3}}\sum _{k=0}^{N}{\frac { \left( -1 \right) ^{k}}{ \left( k+1/2
\mbox{} \right) ^{2\,m+3}}\sum _{j=0}^{p-2}{\frac {\Gamma \left( 2\,j+2 \right)  \left( -1 \right) ^{j}}{ \left( k+1/2 \right) ^{2\,j}{\pi}^{2\,j}}}}  }\,.
\label{S2N}
\end{equation}
The inner sum of \eqref{S2N} can be formally identified in terms of hypergeometric functions as follows

\begin{align} \nonumber
\sum _{j=0}^{p-2}{\frac {\Gamma \left( 2\,j+2 \right)  \left( -1 \right) ^{j}}{ \left( k+1/2 \right) ^{2\,j}{\pi}^{2\,j}}}&={\mbox{$_3$F$_0$}(1,1,3/2;\,\ ;\,\,{\frac {-4}{ {\pi}^{2}\,\left( k+1/2 \right) ^{2}}})}
\mbox{}\\
&+{\frac {\Gamma \left( 2\,p \right)  \left( -1 \right) ^{p}}{ \left( k+1/2 \right) ^{2\,p-2}
\mbox{}{\pi}^{2\,p-2}}{\mbox{$_3$F$_0$}(1,p,1/2+p\,;\,\ ;\,\,{\frac {-4}{ {\pi}^{2}\left( k+1/2 \right) ^{2}}})}}\,.
\label{Sx}
\end{align}

The first $_3F_0(...;.)$ in \eqref{Sx} represents the sum on the left-hand side corresponding to $p\rightarrow \infty$, and the second $_3F_0(...;.)$ represents a correction term for finite values of $p$. Although rigorously divergent when written in series representation, each of these functions can have meaning assigned to it through recourse to Meijer's G-function, and contour integration. Formally, for both $x$ and $q$ complex, (see \cite[Eq. 5.3(1) ]{Luke}) 

\begin{align} \nonumber
\displaystyle  {\mbox{$_3$F$_0$}(1,q,1/2+q\,;\,\ ;\,x)}=&{\frac {{2}^{2q}\,
}{ 2\,\sqrt{\pi}\,\Gamma \left( 2\,q \right) }}G^{1, 3}_{3, 1}\left(-x\, \Big\vert\,^{\displaystyle 0, 1-q, 1/2-q;}_{\displaystyle 0;}\right) \\
=&{\frac {{2}^{2\,q}}{ 2\sqrt{\pi}\,\,\Gamma \left( 2\,q \right) }G^{3, 1}_{1, 3}\left(-\frac{1}{x}\, \Big\vert\,^{\displaystyle 1;}_{\displaystyle 1, 1/2+q, q;}\right)}\,,
\label{Heq}
\end{align}
and we are particularly interested in the case that $x=-4/(\pi^2\,(k+1/2)^2)$ for both general integer values of $q$ as well as $q=1$ in particular. From the theory of G-functions, the final term in \eqref{Heq} can be identified (see \cite[Eq. 5.2(7)]{Luke})


\begin{align} \nonumber
\displaystyle G^{3, 1}_{1, 3}\left(-\frac{1}{x}\, \Big\vert\,^{\displaystyle 1;}_{\displaystyle 1, 1/2+q, q;}\right)=&-{\frac { 8\sqrt{\pi}\,\,\Gamma \left( -2+2\,q \right) }{{2}^{2\,q}\,x}{\mbox{$_1$F$_2$}(1;\,2-q,3/2-q;\,{x}^{-1})}}\\ \nonumber &
+{\pi}^{3/2}\csc \left( \pi\,q \right)  \left( -{x}^{-1} \right) ^{q}{\mbox{$_0$F$_1$}(\ ;\,1/2;\,{x}^{-1})}\\&
-2\,{\pi}^{3/2} \left( -{x}^{-1} \right) ^{1/2+q}\,\sec \left( \pi\,q \right)\,{\mbox{$_0$F$_1$}(\ ;\,3/2;\,{x}^{-1})}\,, 
\label{MeijId}
\end{align}
a result that is valid for $|x|>1$ and complex values of $q$.\newline

{\bf Remark:} The act of generalizing integer $p\rightarrow q$, evaluating the result for general values of $q$ and then reducing the result employing $q\rightarrow p$ as a limit, is known as ``regularization". The functions in \eqref{MeijId} are rooted in confluent hypergeometric series (related to Lommel functions) that possess an infinite radius of convergence. This is significant because the identification $x=-4/(\pi^2\,(k+1/2)^2)$ will eventually result in a series over the index $(k+1/2)^2$.\newline


We now take the limit $q\rightarrow p\,$, and, recognizing that both $k$ and $p$ are integers, each of the terms in \eqref{MeijId} can be simplified; thus \eqref{Heq} becomes
\begin{align} \nonumber
&\displaystyle {\mbox{$_3$F$_0$}(1,p,1/2+p\,;\,\ ;\,{\frac {-4}{{\pi}^{2} \left( k+1/2 \right) ^{2}}})}={\frac {{\pi}^{2} \left( k+1/2 \right) ^{2}
}{\Gamma \left( 2\,p \right) }}\sum _{j=0}^{p-2}\Gamma \left( 2\,p-2\,j-2 \right)  \left( k+1/2 \right) ^{2\,j}{\pi}^{2\,j} \left( -1 \right) ^{j}
\\&
-{\frac {{\pi}^{2\,p} \left( -1 \right) ^{p} \left( k+1/2 \right) ^{2\,p}}{\Gamma \left( 2\,p \right) }\sum _{j=0}^{\infty }{\frac {\psi \left( 2\,j+1 \right)  \left( -1 \right) ^{j} \left( k+1/2 \right) ^{2\,j}{\pi}^{2\,j}}{\Gamma \left( 2\,j
\mbox{}+1 \right) }}}
\mbox{}-{\frac { \left( -1 \right) ^{p+k}{\pi}^{2\,p+1} \left( k+1/2 \right) ^{2\,p}}{2\,\Gamma \left( 2\,p \right) }}\,,
\label{Hp4}
\end{align}

which, with recourse to \eqref{Math1b} can be rewritten as
\begin{align} \nonumber
\displaystyle \frac {\Gamma \left( 2\,p \right) }{{\pi}^{2\,p}}\mbox{$_3$F$_0$}(1&,p,1/2+p;\,\ ;{\frac {-4}{{\pi}^{2} \left( k+1/2 \right) ^{2}}})=
- \left( -1 \right) ^{p} \left( k+1/2 \right) ^{2\,p}\sum _{j=1}^{p-1}{\frac { \left( -1 \right) ^{j}\Gamma \left( 2\,j \right) {\pi}^{-2\,j}}{ \left( k+1/2 \right) ^{2\,j}}}
\mbox{}\\&+ \left( -1 \right) ^{p+k} \left( k+1/2 \right) ^{2\,p}{\rm Si} \left( \pi\,(k+1/2) \right) 
\mbox{}-\pi\, \left( -1 \right) ^{p+k} \left( k+1/2 \right) ^{2\,p}/2
\label{Gans2}
\end{align}

({\bf Remark:} As noted previously, \eqref{Hp4} must vanish as $p\rightarrow \infty$ consistent with the common denominator factor $\Gamma(2p)$.) Since $p>0$ is arbitrary, let $p=1$ in \eqref{Hp4}, to identify

\begin{align} 
\displaystyle{\mbox{$_3$F$_0$}(1,1,3/2;\ ;{\frac {-4}{ {\pi}^{2}\left( k+1/2 \right) ^{2}}})}&={\pi}^{2} \left( k+1/2 \right) ^{2} \left( \sum _{j=0}^{\infty }{\frac {\psi \left( 2\,j+1 \right)  \left( -1 \right) ^{j} \left( k+1/2 \right) ^{2\,j}{\pi}^{2\,j}}{\Gamma
\mbox{} \left( 2\,j+1 \right) }}
+\frac{ \left( -1 \right) ^{k}{\pi}}{2}\right),
\label{Hp1}
\end{align}
a convergent series, which, with the help of \eqref{Math1b} - and true for any choice of integer $p$ - becomes

\begin{equation}
\displaystyle {\mbox{$_3$F$_0$}(1,1,3/2;\,\ ;\,{\frac {-4}{ \left( k+1/2 \right) ^{2}{\pi}^{2}}})}=-{\pi}^{2} \left( -1 \right) ^{k} \left( k+1/2 \right) ^{2}
\mbox{} \left( {\rm Si} \left( \pi\, \left( k+1/2 \right)  \right) -\pi/2 \right) \,.
\label{Hp2}
\end{equation}

({\bf Remark:} Inserting \eqref{Hp4} and \eqref{Hp2} into \eqref{Sx} yields a trivial identity, reinforcing confidence that the somewhat unusual analysis methods employed thus far are valid.)\newline

Comparing \eqref{Gans2} and \eqref{Gans4} now gives
\begin{equation}
\displaystyle {E_{2p}} \left( -i\pi\, \left( k+1/2 \right)  \right) +{ E_{2p}} \left( i\pi\, \left( k+1/2 \right)  \right) =-
\mbox{}\,\frac { 2\left( -1 \right) ^{k}}{\pi\, \left( k+1/2 \right) }{\mbox{$_3$F$_0$}(1,p,1/2+p;\,\ ;\,{\frac {-4}{{\pi}^{2} \left( k+1/2 \right) ^{2}}})}\,,
\label{Conj1}
\end{equation}

and with reference to \eqref{Heq} and \cite[Eq. (5.2(1)] {Luke}, the right-hand side of \eqref{Conj1} can be rewritten as a well-defined contour integral:


\begin{align} \nonumber
\displaystyle \mbox{$_3$F$_0$}&(1,p,1/2+p\,;\,\ ;\,{\frac {-4}{{\pi}^{2} \left( k+1/2 \right) ^{2}}})=\\
&-i\frac {{2}^{2\,p-2}
}{{\pi}^{3/2}\Gamma \left( 2\,p \right) }\int_{c-i\infty }^{c+i\infty }\!\Gamma \left( 1-t \right) \Gamma \left( 1/2+p-t \right) 
\mbox{}\Gamma \left( p-t \right) \Gamma \left( t \right)  \left( \pi/2 \right) ^{2\,t} \left( k+1/2 \right) ^{2\,t}\,{\rm d}t
\label{Thus}
\end{align}
where $0<c<1$. \newline

Keeping \eqref{Cp} in mind, it is convenient to translate the contour of integration in \eqref{Thus} such that $m<c<m+1$, (denoted by $c_m$) and, with $p>m$, by adding the appropriate residues we obtain

\begin{align} \nonumber
\displaystyle \frac {{E_{2p}} \left( -i\pi \left( k+1/2 \right)  \right) +{ E_{2p}} \left( i\pi \left( k+1/2 \right)  \right) }{ \left( k+1/2
\mbox{} \right) ^{2\,p+2\,m}}&={\frac {i \left( -1 \right) ^{k}}{\pi\,\Gamma \left( 2\,p \right) }\int_{c_m-i\infty }^{c_m+i\infty }\!{\frac {\Gamma \left( 2p-2t \right) {\pi}^{2t} 
}{\sin \left( \pi\,t \right)\left( k+1/2 \right) ^{-2\,t+2\,p+2\,m+1} }}\,{\rm d}t}\\
&\hspace{-30pt}+\frac {2 \left( -1 \right) ^{k} 
\mbox{}}{\pi\,\Gamma \left( 2\,p \right) \,\left( k+1/2 \right) ^{2p+2\,m-2\,j+1} }\sum _{j=1}^{m} \left( -1 \right) ^{j}{\pi}^{2\,j}\Gamma \left( 2\,p-2\,j \right)\,.
\label{Q1a}
\end{align}

Now, employing \eqref{Fsum1a} and \eqref{BetaDef}, substitute \eqref{Q1a} into \eqref{Cp}, split one of the sums to cancel, and reverse a remaining sum, giving
\begin{align} \nonumber
\displaystyle  &\left( {2}^{2\,m+1}-1 \right) \zeta \left( 2\,m+1 \right) ={\frac {i \left( -1 \right) ^{p}}{\pi^{2p+1}}\int_{c_m-i\infty }^{c_m+i\infty }\!{\frac {\Gamma \left( 2\,p-2\,v \right) {\pi}^{2\,v}
}{\sin \left( \pi\,v \right) }}\sum _{k=0}^{\infty } \left( -1 \right) ^{k} \left( k+1/2 \right) ^{2\,v-2\,p-2\,m-1}\,{\rm d}v}\\
&- \left( -1 \right) ^{m}{\pi}^{2\,m}\sum _{j=0}^{p-2-m}{\frac {\Gamma \left( 2\,j+2 \right) \mathcal{E} \left( 2\,j+2+2\,m \right) }{\Gamma \left( 2\,j+3+2\,m \right) }}
\mbox{}- \left( -1 \right) ^{m}{\pi}^{2\,m}\sum _{j=0}^{m}{\frac {\mathcal{E} \left( 2\,j \right) \psi \left( 1-2\,j+2\,m \right) }{\Gamma \left( 1-2\,j+2\,m \right) \Gamma \left( 2\,j+1 \right) }}\,.
\label{Cpa4}
\end{align}


As noted earlier, the sum of the first two terms on the right-hand side of \eqref{Cpa4} is independent of $p$. One possibility now is to choose $c_m=m+1/2$ along with $p=m+1$. With an obvious change of integration variables, the final result is given by \eqref{Cp5a}.

\end{appendices}
\end{flushleft}
\end{document}